\documentclass[reqno, 10pt]{amsart}

\usepackage[english]{babel}
\usepackage{enumitem}
\usepackage{graphicx}
\usepackage{graphicx}
\usepackage{color}
\usepackage[dvipsnames]{xcolor}
\usepackage{hyperref}
\usepackage{mathtools}
\usepackage{comment}
\usepackage{amsmath}
\usepackage{amssymb}
\usepackage{amsfonts}
\usepackage{amsthm}
\usepackage{mathrsfs}
\usepackage{soul}
\usepackage{float}

\newcommand{\N}{\mathbb{N}}
\newcommand{\R}{\mathbb{R}}

\newcommand{\Z}{\mathbb{Z}}

\newcommand{\f}{\rightarrow}

\newcommand{\JSJ}{\textit{JSJ}}

\newcommand{\fix}{\mathsf{Fix}}

\newcommand{\ent}{\operatorname{Ent}}

\newcommand{\Ent}{\operatorname{Ent}}

\newcommand{\diam}{\operatorname{diam}}
\newcommand{\Stab}{\operatorname{Stab}}

\newcommand{\Ric}{\operatorname{Ricci}}

\newcommand{\ax}{\mathsf{Ax}}

\newtheorem{thm}{Theorem}[section]
\newtheorem{thmintro}{Theorem}

\newtheorem{prop}[thm]{Proposition}
\newtheorem{lem}[thm]{Lemma}
\newtheorem*{lem*B}{Lemma B}
\newtheorem*{lem*C}{Lemma C}
\newtheorem*{lem*D}{Lemma D}
\newtheorem*{lem*E}{Lemma E}
\newtheorem*{prop*A}{Proposition A}

\newtheorem{corintro}[thmintro]{Corollary}

\theoremstyle{remark}

\newtheorem{rmk}[thm]{Remark}
\newtheorem*{rmkintro}{Remark}
\newtheorem*{rmk*}{Remark}

\theoremstyle{definition}
\newtheorem{defn}[thm]{Definition}
\newtheorem*{defn*}{Definition}

\begin{document}

\title[Entropy and finiteness of groups with acylindrical splittings]{Entropy  and finiteness of groups \\with acylindrical splittings}
\author[F.Cerocchi]{Filippo Cerocchi}
\date{\today}
\address{F. Cerocchi, Max-Planck-Institut f\"ur Mathematik, Vivatsgasse 7, 5311, Bonn.}
\author[A. Sambusetti]{Andrea Sambusetti}
\address{A. Sambusetti, Dipartimento di Matematica \lq\lq G. Castelnuovo'', Sapienza Universit\`a di Roma, P.le Aldo Moro 5, 00185, Roma.\newline}
\email{\newline F. Cerocchi: fcerocchi@mpim-bonn.mpg.de\, - \, fcerocchi@gmail.com \newline
A. Sambusetti: sambuset@mat.uniroma1.it}

\maketitle

\vspace{-3mm}
{\bf Abstract.}
We prove that  there exists a positive, explicit function $F(k, E)$  such that,  for any  group $G$ admitting a  $k$-acylindrical splitting and any generating set $S$ of $G$ with $\Ent(G,S)<E$, we have $|S| \leq F(k, E)$. 
We deduce corresponding finiteness results for  classes of groups possessing acylindrical splittings and acting geometrically with bounded entropy: for instance,  $D$-quasiconvex $k$-malnormal amalgamated products acting on $\delta$-hyperbolic spaces or on $CAT(0)$-spaces  with entropy bounded by $E$. 
A number of finiteness results for interesting families of Riemannian or metric spaces with bounded entropy and diameter also follow:  Riemannian 2-orbifolds,  non-geometric $3$-manifolds, higher dimensional graph manifolds and cusp-decomposable manifolds,   ramified coverings and, more generally,  CAT(0)-groups with negatively curved splittings.

\tableofcontents


\vspace{-10mm}
\section*{Introduction}

 In this paper  we are interested in finitely generated groups $G$ {\em admitting $k$-acylindrical splittings}, that is isomorphic to the fundamental group of a graph of groups  such that
  the action of $G$ on the corresponding Bass-Serre tree is 
 (non elementary and) 
  $k$-acylindrical.
 The notion of  acylindricity  is due to Sela in \cite{Sela}, and arises naturally in the context of Bass-Serre theory. 
 It is a geometric translation of   the notion of {\em malnormal amalgamated product}, introduced by Karras and Solitar \cite{karso} at the beginning of the seventies (see Section \S1 for details). \linebreak
 We recall that an action without inversions of a finitely generated group on a simplicial tree $\mathcal T$ is said to be {\it $k$-acylindrical} if the fixed point set of any element $g\in G$ has diameter at most  $k$.
 
%
The more elementary examples of groups possessing a $k$-acylindrical splitting  are  the free products
and the fundamental groups of compact surfaces of negative Euler characteristic, but the class is considerably larger and encompasses several interesting classes of amalgamated groups  which naturally arise in Riemannian and metric geometry, as we shall see.
During the years, the existence of acylindrical actions on simplicial trees has been mainly used to prove some accessibility results (\cite{Sela}, \cite{kaweid2}).  More recently, there has  been an increasing interest on groups acting acylindrically on Gromov hyperbolic spaces (see \cite{osin}, \cite{MiOs}, \cite{sisto}  and references therein for this more general notion of acilindricity). 
We shall instead focus on {\em growth} and {\em finiteness} properties of such groups and of some classes of spaces arising as quotients of geometric actions of these groups. 
\vspace{1mm}

Recall that the \textit{entropy}, or exponential growth rate,   of a group $G$  with respect to a finite generating set  $S$ is defined as 
$$\ent(G,S):=\lim_{n\f+\infty}\frac{1}{n} \log  \left|  S^n \right|$$
where $|S^n|$ is the cardinality of the ball of radius $n$ centered at the identity element, with respect to the word metric   relative to $S$. We will also deal with   groups $G$ acting discretely by isometries on general (non-compact)  proper metric spaces $Y$; in this case, the entropy of the action is defined as 
$$\ent(G \curvearrowright Y):=\limsup_{R\f+\infty}\frac{1}{R} \log  \left|   B_Y(y_0,R) \cap Gy_0 \right|$$
where $B_Y(y_0,R)$ denotes the ball  of radius $R$ centred at $y_0$
(the limit is clearly independent from the choice of the base point $y_0 \in Y$). 
When $X$ is a closed Riemannian manifold it is customary to call {\em (volume-)entropy} of $X$ the entropy of $G=\pi_1(X)$ acting on its Riemannian universal covering space $Y=\widetilde X$; this number coincides with the exponential growth rate of the volume of balls in $X$, and  it is well known that it equals, in non-positive curvature,  the {\em topological entropy} of the geodesic flow on the unitary tangent bundle of $X$, cp. \cite{manning}.    
We  extend this terminology to any compact metric space $X$ obtained as the quotient $X=G \backslash Y$ of a simply connected, geodesic space $Y$ by a discrete group of isometries (possibly, with fixed points). 
%
We will come back shortly to  the analytic information encrypted in this asymptotic invariant for Riemannian manifolds.
\vspace{1mm}

 The first main result of this paper, from a group-theoretic point of view, is the following:
 
 \begin{thmintro}[Entropy-Cardinality inequality]\label{ent-card_ineq}
  Let $G$ be a group acting by automorphisms  without inversions,  non-elementarily and  $k$-acylindrically on a   tree.\linebreak
 For any finite generating set $S$ of $G$ of cardinality $n$, we have 
 \vspace{-2mm}
 
 \begin{equation}\label{eqEC}
 \ent(G,S) \ge \frac{\log \left(\sqrt[32]{n}-1\right)}{80(k+3)}
 \end{equation}
 \end{thmintro}

 \begin{rmkintro} The above inequality is far from being optimal, and makes sense only for $n \gg 0$. The meaning is qualitative: the entropy $\ent(G,S)$ diverges as $|S|$ becomes larger and larger.  
\end{rmkintro}
A first remarkable example of  Entropy-Cardinality inequality was given by \linebreak Arzhantseva-Lysenok  \cite{arly}:
  for any given  hyperbolic group $G$ there exists a
  constant $\alpha(G)$
  such that for any non-elementary, finitely generated subgroup
   $H$ and any finite generating set $S$ of $H$ the  inequality   $\ent(H,S)\ge \log\left(\alpha(G)\cdot |S|\right)$ holds  \footnote{It is worth noticing here that given a finitely generated group $G$ acting $k$-acylindrically on a simplicial tree $\mathcal T$, any finitely generated subgroup $H<G$ also acts $k$-acylindrically on  the minimal subtree $\mathcal T_H$, globally preserved by $H$; therefore Theorem \ref{ent-card_ineq} also holds for all subgroups of $G$, provided that their action  on $\mathcal T $ is still non-elementary.}. 
  One interest of similar inequalities  is that they generally represent  a step forward to prove (or disprove) the realizability of the {\em algebraic entropy}
  \footnote{The algebraic entropy $\Ent_{alg} (G)$ of a finitely generated group $G$ is  defined as the infimum of the entropies $\Ent(G,S)$, when $S$ varies among all possible, finite generating sets for $G$.}  for a group $G$. 
 On the other hand, the Entropy-Cardinality inequality proved in this paper has a different theoretical meaning, since the   cardinality of $S$ is bounded in terms of a {\em universal} function, only depending on the entropy  and   the acylindricity constant $k$, and not on the  group $G$ itself.
  
The idea of proof of (\ref{eqEC}) is relatively elementary  and based on the construction of free subgroups of $G$ of large rank.
Any collection of hyperbolic elements of $G$ admitting  disjoint  domains of attraction generates a free subgroup  (namely, a \textit{Schottky subgroup})  by a classical ping-pong argument. So, the strategy for Theorem \ref{ent-card_ineq} is to show that from any  sufficiently large generating set $S$ one can produce  large collections (compared to $|S|$)  of hyperbolic elements  with suitable configurations of axes and uniformly bounded $S$-lengths. 
However, quantifying this idea turns out to be a rather complicate combinatorial problem,  since we need to  control  the mutual positions of the axes  as well as the translational lengths  of all the hyperbolic elements under consideration. 
The first part of the paper (\S\ref{subsecbasic}-\S\ref{ssextraction}), is entirely devoted to developing  the combinatorial tools needed to prove this  inequality.

\vspace{1mm}
  In the second part of the paper we focus on algebraic and geometric applications of the Entropy-Cardinality inequality.  As an immediate consequence,   we get general finiteness results for abstract groups with $k$-acylindrical splittings and uniformily bounded entropy, provided we know that they possess a complete set of relators of uniformily bounded length (cp. Theorem  \ref{finiteness_short_relators}  in Section \S\ref{ssapplications}). The following are particularly interesting cases:

\begin{corintro}\label{deltahyp} The number of  isomorphism classes of marked, $\delta$-hyperbolic groups $(G,S)$ with a non-elementary $k$-acylindrical splitting and  satisfying  $\ent(G,S) \leq E$ is finite, bounded by an explicit function $M(k, \delta, E)$.
\end{corintro}

\begin{corintro}\label{quasiconvex} The number of  groups $G$ admitting a non-elementary $k$-acylindrical splitting, with a  $D$-quasiconvex action on: \\
\noindent (i) either  some (proper, geodesic)  $\delta$-hyperbolic space $(X,d)$,\\
\noindent (ii) or on some CAT$(0)$-space $(X,d)$,\\
and satisfying    $\Ent(G \curvearrowright X) \leq E$  is finite. Their number is  bounded by a  function of $k, \delta, D, E$  in case (i), and of  $k, D, E$  in case (ii).
\end{corintro}

\noindent (We stress the fact that, in the above corollary, the hyperbolic or $CAT(0)$-spaces $X$ the group $G$ acts on are not supposed to be fixed).
  \vspace{1mm}

A  typical case where   a group $G$ admits an action on a $CAT(0)$-space and a $k$-acylindrical splitting occurs  when $G$ is the fundamental group of a space \linebreak $X=X_1 \sqcup_{Z} X_2$  which is the gluing of   two locally $CAT(0)$-spaces $X_i$ along two  isometric, locally convex  subspaces $Z_i \cong Z$
 (or  $X=X_0 \sqcup_{\phi} $ is obtained by identifying two such  subspaces $Z_i \subset X_0$ to each other by an isometry $\phi$), and the resulting space $X$ is locally,  negatively curved around $Z$.
Namely, we will say  that $X$ has a {\em negatively curved splitting} if   the subspace $Z$ obtained by identifying $Z_1$ to $Z_2$ has a neighbourhood $U(Z)$ in $X$  such that $U(Z)\setminus Z$ is  a locally $CAT(-\kappa)$-space, for some $\kappa >0$. 
The fact that $Z$ possesses such a neighborhood ensures
\footnote{Notice that some form of {\em strictly} negative curvature must be assumed to deduce that   $\pi_1(Z)$ is  malnormal in  $\pi_1(X)$, and it is not sufficient to require that $Z$ is a locally convex subset of a locally $CAT(0)$-space $X$ of rank $1$. 
A counterexample is provided, for instance, by a $3$-dimensional  irreducible manifold $X$ with non trivial JSJ splitting and  one  component of hyperbolic type:  by  \cite{leeb},  $X$ can be given a  non-positively curved metric  of rank $1$ in the sense of \cite{Bal},  and  possesses a totally geodesic, embedded torus $Z$ whose fundamental group is not malnormal in $\pi_1(X)$.}   
that $\pi_1(Z)$ is a malnormal subgroup in each $\pi_1(X_i)$, and therefore $\pi_1(X)$ has a 1-acylindrical splitting: we refer to   Section \ref{ssnegcurvsplittings} and to the Appendix \ref{appcat(0)} for details.
We then have:

\begin{thmintro} \label{thmcat(0)}
 The number of homotopy types of compact, locally  $CAT(0)$-spaces $X$ admitting a non-trivial negatively curved splitting, satisfying $\Ent(X)<E$ and $\diam(X)<D$ is finite. 
\end{thmintro}


 \begin{corintro} \label{cormanifoldsk<0}
 There exist only finitely many non-diffeomorphic closed,  non-positively  curved manifolds $X$ of dimension different from 4 admitting a non-trivial  negatively curved  splitting and satisfying $\Ent(X)<E$ and $\diam(X)<D$. In dimension 4, the same is true up to homotopy equivalence.
\end{corintro}

\noindent It is worth noticing that Corollary \ref{cormanifoldsk<0}  holds more generally for manifolds $X$ with metrics of curvature of {\em any} possible sign, provided that $X$ also admits a non-trivial, negatively curved splitting (this follows from the combination of Theorem \ref{thmcat(0)} with Lemma \ref{SPorbifold}, as in the proof of Theorem \ref{fntn_orb} for non-triangular $2$-orbifolds).


\vspace{1.5mm}
 
 
 The question whether a given family of Riemannian manifolds contains only a finite number of topological types has a   long history: the ancestor of all finiteness results is probably Weinstein's theorem   \cite{wein}  
 	on  finiteness of the homotopy types of pinched, positively curved, even dimensional manifolds.
 A  few years later,  Cheeger's celebrated Finiteness Theorem appeared, for closed Riemannian manifolds with bounded   sectional curvature and, respectively, lower and upper bounds on volume and diameter \cite{chee},  \cite{pet}.
Several generalizations  \cite{GP}, \cite{GPW}
 with relaxed  assumptions on the curvature
 then followed,  until, in the nineties, the attention of geometers turned to Riemannian manifolds satisfying a lower bound on the Ricci curvature, driven by  Gromov's Precompactness Theorem. 
 Substantial progresses in  understanding  the diffeomorphism type of Gromov-Hausdorff limits  and  the local structure of manifolds under  lower Ricci curvature bounds  were then made  --by no means trying to be exhaustive--    by 
Anderson-Cheeger \cite{anderson-cheeger}, Cheeger-Colding \cite{ch-co}  and, more recently by 
Kapovich-Wilking \cite{ka-wi}  (see also  Breuillard-Greene-Tao's  work  \cite{BGT},   for a more group-theoretical  approach to the generalized Margulis' Lemma under packing conditions --a macroscopic  translation of a lower Ricci curvature bound).

 \pagebreak
 
 Corollary \ref{cormanifoldsk<0} represent an attempt to get rid of lower curvature bounds, at least in non-positive curvature, replacing it only by a bound of an asymptotic invariant.  
 Recall that, for  a closed Riemannian manifold  $X$,  a lower  bound  of the Ricci curvature $\Ric_X \geq -(n-1)K^2$ implies  a corresponding upper bound of the entropy $\Ent(X) \leq (n-1)K$,  by the classical volume-comparison theorems of Riemannian geometry. However the entropy,  being an asymptotic invariant, only depends on the large-scale geometry of the universal covering $\widetilde X$, and can  be seen as an averaged version of the curvature (this can be given a precise formulation in negative curvature by  integrating  the Ricci curvature on the unitary tangent bundle of $X$ with respect to a suitable measure, cp. \cite{knieper}).
 Therefore, the condition $\Ent(X)<E$  is much weaker than a  lower bound on the Ricci curvature.  
To get a glimpse of the difference, remark that the class of closed, Riemannian  manifolds (of dimension $n\ge 3$) with uniformily bounded entropy and diameter is not precompact with respect to the Gromov-Hausdorff distance (see \cite{rev}, Remark 2),  and neither is the  family of Riemannian  structures with  uniformily bounded entropy and diameter on any given $n$-dimensional manifold $X$ (see \cite{revthese}, Example 2.29).

The first  results  about families of Riemannian metric and metric-measured spaces with uniformily bounded entropy, such as lemmas \`a-la-Margulis, finiteness and compactness results etc., were given by Besson-Courtois-Gallot's in \cite{BCG3} (yet unpublished), and are the object of  \cite{BCGS}.
Other {\em local} topological rigidity results under entropy bounds  have  recently appeared in the authors' \cite{CS}.

Under this perspective,   Corollary  \ref{cormanifoldsk<0}  might be compared with the well-know 
finiteness result for negatively  curved $n$-manifolds $X$  with uniformily bounded diameter and  {\em sectional curvature}  $ K(X) \geq -k^2$ 
(which follows from a version of the  Margulis' Lemma in non-positive curvature, as stated for instance in \cite{buka}, and from the aforementioned Cheeger's  finiteness theorem).
It is a challenging open question to know whether  the conclusion of Corollary \ref{cormanifoldsk<0}  extends to  {\em all}  closed, negatively curved manifolds with uniformily bounded entropy and diameter.

 	
 	
 	
\vspace{2mm}


 The Entropy-Cardinality inequality becomes a powerful tool, when applied to families of Riemannian manifolds enjoying strong topological-rigidity properties.  \linebreak To illustrate this fact, we present now some basic examples of application  of  Theorem \ref{ent-card_ineq} to particular  classes of spaces whose groups  naturally possess acylindrical splittings and presentations with an uniform bound on the acylindricity constant and  on the length of relators.  
\vspace{3mm}

{\bf  A. Two-dimensional orbifolds of negative orbifold  characteristic.}

 \noindent {\em Orbifolds} were introduced by  Satake \cite{sat} in the late fifties under the name of {\em $V$-manifolds} and later popularized by Thurston  (\cite{Thurst}) who used them to show the existence of locally homogeneous metrics on Seifert fibered manifolds. Generally speaking, $n$-dimensional orbifolds 
 are mild generalization of manifolds, whose points have neighborhoods  modeled on the quotient of $\R^n$  (or on the upper half space $\R^n_+$) by the action of a finite group of transformations; we refer to section \S\ref{ss2dimorbifolds} for precise definitions and  isomorphisms of 2-orbifolds. 
In the '80s Fukaya introduced the equivariant Gromov-Hausdorff distance (\cite{fuk}), and used it to study Riemannian orbifolds. Since then, several authors gave attention to spectral and finiteness results on Riemannian orbifolds (see for instance \cite{bor}, \cite{farsi}, \cite{stan}, \cite{prostan}, \cite{pro}), possibly because of their application to string theory (\cite{ALR}). 
We show:

\begin{thmintro}\label{fntn_orb}${}$
Let $\mathscr O^2_-(E,D)$ be the class of Riemannian, compact, $2$-orbifolds (with or withour boundary) with conical singularities and $\chi_{orb}(\mathcal O)\le0$, satisfying \linebreak $\ent(\mathcal O)\le E$ and $\diam(\mathcal O)\le D$. This class contains only  a finite number of  isomorphism types.
\end{thmintro}

 \vspace{-1mm}

 \noindent We stress  the fact that the orbifold metrics in the class $\mathscr O^2_-(E,D)$  under consideration  {\em are not} supposed to be negatively or nonpositively  curved.

 \noindent Notice that the analogous result for compact {\em surfaces}  easily follows  from  basic estimates of the algebraic entropy of a surface group,  together with the aforementioned  Gromov's inequality 
 \small $$\ent(X)\diam(X)\ge \frac{1}{2}\ent_{alg}(\pi_1(X))$$
 \normalsize
Actually, it is well known that the algebraic entropy of a compact  surface $X$ of genus $g$ with $h$ boundary component  is bounded from below by $\log (4g+3h-3)$  (cp. \cite{delH}), therefore  $\ent(X)$ and $\diam(X)$ bound   $g$ and $h$.
The orbifold case is significantly more tricky:  we use the fact that {\em non-triangular} orbifolds of negative Euler characteristic always admit  a $2$-acylindrical splitting 
 (a proof of this is given in the Appendix \ref{apporbifold}, Proposition \ref{splitting2orbifolds}), so we can apply  the Entropy-Cardinality inequality to particular, well-behaved presentations of the orbifold  groups.
On the other hand, triangular orbifolds do not admit such splittings, and we are forced to  a direct  computation, using arguments from classical small cancellation theory.

\begin{rmkintro} The above finiteness theorem marks a substantial difference with the analogue question in geometric group theory: in fact, the number of  $2$-orbifold groups $G$ admitting a generating set $S$ such that $\Ent(G,S) \leq E$ {\em is not} finite (at least, without any additional, uniform hyperbolicity assumption on the groups $G$). 
Actually, on {\em any} topological surface $S$ of genus $g$ with $k$ conical points of orders $p_1,...,p_k$ there always exists a generating set $S$ of cardinality at most  $2g+k$, such that $\Ent(G,S)$ is smaller than the entropy of the free group on $S$, independently from the choice of the orders $p_1,...,p_k$ (we thank R. Coulon for pointing out this fact to us).
The reason for this difference is that, on Riemannian orbifolds,  any torsion element $g \in G$ has a fixed point on $\widetilde X$, which gives rise to arbitrarily small loops increasing substantially the entropy; this does not happen for the action of $G$ on its (non-simply connected) Cayley graph ${\mathcal G} (G,S)$. Notice also that the existence of torsion elements with  unbounded orders prevents $(G,S)$ to be $\delta$-hyperbolic, for any fixed $\delta$. 
 \end{rmkintro}
\vspace{3mm}

{\bf B. Non-geometric $3$-manifolds.}\\
A compact  $3$-manifold, possibly with boundary, is called {\em non-geometric}
\footnote{We use here  the term ``geometric'' as in the original definition given in \cite{thur82}; in the  case of manifolds with boundary, variations on this definition are possible and suitable for other purposes (i.e. uniqueness of the model geometries on each piece),  see for instance \cite{Bon}.}
if its interior cannot be endowed with a complete Riemannian metric locally isometric to one of the eight $3$-dimensional  complete, maximal, homogeneous model geometries:  $ {\mathbb E}^3, {\mathbb S}^3, {\mathbb H}^3, {\mathbb S}^2 \! \times \! {\mathbb R}, {\mathbb H}^2 \! \times \! {\mathbb R},  \mathbb H^2\widetilde\times\mathbb R, Nil$ and $Sol$.  
We restrict our attention to  orientable  manifolds,  for the sake of simplicity
\footnote{The geometrization conjecture for non-orientable $3$-manifolds with boundary being not yet clearly estblished, as far as the authors know.},
{\em without spherical boundary components} (since, clearly, punctures cannot be detected by the fundamental group in dimension 3).\\
By classical results of $3$-dimensional topology and by the solution of the Geometrization Conjecture, any such manifold  $X$ falls into one of the following mutually disjoint classes (as explained, for instance, in \cite{AFW}):

 (i) either $X$ is {\em not prime}, and is different from $P^3{\mathbb R} \# P^3{\mathbb R}$ (the only compact, non-prime manifold without spherical boundary components admitting a geometric structure, modelled on ${\mathbb S}^2 \! \times \! {\mathbb R}$);
 
(ii) or $X$ is {\em irreducible}, has a non-trivial $JSJ$ decomposition and is not  finitely covered by a torus bundle (in which case, it would admit a $Sol$-geometry). \\
In the first case, the fundamental group of $X$ admits a non-elementary, $1$-acylindrical splitting corresponding to the prime decomposition, while in case (ii) the JSJ-decomposition induces a (at most) $4$-acylindrical splitting of $\pi_1(X)$, by  \cite{WiZa}
  (cp. also Section \S4 in  \cite{CS}, and \cite{cer_pi1} for details about the degree of acylindricity of the splitting over  the abelian subgroups corresponding to the JSJ tori, according to the different types of adjacent JSJ-components).

In \cite{CS} the authors exhamined the  local rigidity  properties of Riemannian, non-geometric $3$-manifolds with {\em torsionless} fundamental group, under  uniformly bounded entropy and diameter.  
Here we consider, more generally, the class  $\mathscr M_{ng}^{\partial}(E,D)$ of  compact, orientable  non-geometric $3$-manifolds  (with possibly empty, non-spherical boundary), possibly with torsion, endowed with Riemannian metrics with entropy and diameter bounded from above by two positive constants $E$ and $D$. 
Acylindricity of the splitting of their fundamental groups is the key to the following:
	\begin{thmintro}\label{glob_fnt}
		The number of isomorphism classes of fundamental groups of  manifolds in  $\mathscr M_{ng}^{\partial}(E,D)$ is  less than $\left(e^{1120\,E\,D} +1\right)^{32}$.
	\end{thmintro}
		
 
Moreover,  the homotopy type (and, in turns, the diffeomorphism type) of  compact $3$-manifolds  without  spherical boundary components   is determined by their fundamental group up to a finite number of choices, by  Johannson and Swarup works \cite{jo2},\cite{swa}  (see the discussion in Section \S\ref{sec3manifolds} for details,  and Theorem \ref{pi1_determines_homeo}  in Appendix \ref{appJS}). We therefore obtain:
 	
	\begin{corintro}
		\label{fnt_HT}
		The number of diffeomorphism types in $\mathscr M_{ng}^{\partial}(E,D)$ is finite.
	\end{corintro}

	\noindent  Notice that, while Theorem \ref{glob_fnt} gives an explicit (albeit ridicolously huge) estimate of the number of groups in $\mathscr M_{ng}^{\partial}(E,D)$, Corollary \ref{fnt_HT} does not  provide any  explicit estimate of the number of diffeomorphism types.  

\begin{rmkintro}
	The bound  we found in Theorem  \ref{glob_fnt} is explicit, but   far from being optimal. Once finiteness for this class is known, one might  try to use more efficient, computer-assisted algorithms to find reasonable estimates of their number. 
\end{rmkintro}

 \vspace{3mm}

{\bf  C. Ramified coverings of hyperbolic manifolds}\\
Another interesting  class of spaces whose fundamental groups admit acylindrical splittings is the one  of cyclic ramified coverings of hyperbolic manifolds.  The construction  is due to  Gromov-Thurston (\cite{GT}) and  represents an important source of examples of  manifolds  admitting pinched, negatively curved metrics  but  not hyperbolic ones. 
 A degree $k$ ramified cover $X_k$ of a hyperbolic manifold $X$ is obtained by excising a  totally geodesic hypersurface with boundary $Z$ in $X$, and then  glueing  several copies of $X \!-\!Z$  along a  ``$k$-paged open book'', whose leaves are copies $Z_i$ of $Z$ joined together at the ramification locus $R=\partial Z$. Any such covering admits a singular, locally $CAT(-1)$ metric, and its fundamental group splits as a free product  of CAT(-1) groups amalgamated  over the fundamental group of the locally convex subspace $Z_1 \cup Z_k$ of $X_k$ (given by two  pages of the book). We will recall in  Section \S2.5  this construction in more detail, and show that these manifolds naturally fall in the class of spaces with negatively curced splittings.
However, we will be  interested in metrics with curvature of {\em any possible sign} on such manifolds. \linebreak 
Namely,  let $ \mathscr R^4 (E,D)$ be the space of all 4-dimensional Riemannian, cyclic ramified coverings  of compact, orientable hyperbolic manifolds, whose entropy and diameter are respectively bounded by $E$ and $D$;  and let  $\mathscr R^{\neq 4}  (E,D)$ the corresponding space   of ramified coverings in dimensions different  from $4$.  
Then: 
 
\begin{corintro}\label{corram}The class  $\mathscr R^4 (E,D)$  contains only finitely many different homotopy types, and   $\mathscr R^{\neq 4} (E,D)$ only  finitely many diffeomorphism  types.
\end{corintro}

\noindent We stress the fact that the manifolds under consideration are (genuine) ramified coverings    of {\em any} possible hyperbolic manifold, and not just of {\em one} fixed hyperbolic manifold. 
\vspace{3mm}

{\bf  D. High  dimensional graph and cusp decomposable manifolds.}

\noindent  High  dimensional graph manifolds have been introduced by Frigerio, Lafont and Sisto in \cite{FLS},  and cusp-decomposable manifolds by Nguyen Phan in \cite{nguyen}. \\
Roughly speaking, an (extended) {\em $n$-dimensional graph manifold} $X$ is obtained gluing together, via affine diffeomorphisms  of their boundaries, several elementary, building blocks $X_i$ which are diffeomorphic to products  $H^{k_i} \times T^{n-k_i}$, where $T^{n-k_i}$ is a $(n-k_i)$-dimensional torus (representing the local {\em fibers} of the graph manifold), and $H^{k_i}$ is a manifold of  dimension $k_i \geq 2$ with toroidal boundary components, obtained from a hyperbolic manifold with cusps by truncating some cusps along (the quotient of)   flat horospheres; 
a block of the form  $H^{2} \times T^{n-2}$ is also called a {\em surface piece}
\footnote{Notice that the above class  does no cover the case of non-geometric $3$-manifolds, since  surface pieces are trivial products, so the blocks do not take into account Seifert fibrations.}.
 A high  dimensional graph manifold $X$ is called {\em irreducible} if none of the fibers of two adjacent blocks represent the same element in $\pi_1(X)$,  and {\em purely hyperbolic} if there are no fibers at all (that is  every piece $X_i$ is a truncated $n$-dimensional hyperbolic manifold). \\
{\em Cusp decomposable manifolds} are defined similarly to purely hyperbolic high  dimensional graph manifolds, but starting from building blocks which are  negatively curved, locally symmetric manifolds with cusps, and glueing, always via  affine  diffeomorphisms\footnote{An {\em affine diffeomorphism}, in this context, is the composition of an isomorphism of the two nilpotent Lie groups composed with the left multiplication by an element of the group. In \cite{nguyen} the author points out  the necessity of realizing the gluings via affine diffeomorphisms to have strong differential rigidity: in fact, Aravinda-Farrell show   in \cite{arfa}  the existence of non-affine gluing maps, for the double of a hyperbolic cusped manifold $X$,  giving rise to a manifold which  is not diffeomorphic to the one obtained by gluing the two copies via the identity map of $\partial X$.}, the boundary  infra-nilmanifolds obtained by truncating the cusps.  
 We say that a $n$-dimensional graph manifold, or a cusp decomposable manifold $X$, is {\em non-elementary} if it is obtained by identifying at least two boundaries (of one or more  building blocks).

High  dimensional graph enjoy strong topological  rigidity properties, as they are aspherical spaces and satisfy the Borel Conjecture in dimension $n\geq 6$ (cf. \cite{FLS}, Thm.3.1 and \S3.4, Rmk.3.7); moreover, the  diffeomorphism type  of (closed) cusp decomposable manifolds, or of all high dimensional graph manifolds whose boundary components do not  belong to surface pieces,  is determined by the fundamental group within the respective classes (cp. \cite{FLS} Thm.0.7 and \cite{nguyen}). 
The importance of considering {\em affine gluings}, in order to get rigidity, was pointed out   by Aravinda-Farrell in \cite{arfa}, where they considered the so called {\em twisted doubles}  ---which also are included in the class defined in \cite{FLS}.

 Non-elementary, irreducible high  dimensional graph manifolds groups possess $2$-acylindrical splittings, while  cusp decomposable manifolds (or purely hyperbolic, high  dimensional graph manifolds) groups have $1$-acylindrical splittings. Therefore, we obtain:

\begin{corintro}
\label{corgraphcusp} $\!\!\!$ Let $\mathscr  G ^{\partial}(E,D)$  (resp. $\mathscr G (E,D)$) be the class of  compact (resp. closed) \linebreak Riemannian,  non-elementary irreducible   high dimensional graph manifolds  with entropy and diameter bounded by $E$, $D$,  and let  $\mathscr C (E,D)$ be the class of closed Riemannian, non-elementary cusp decomposable manifolds satisfying the same bounds:\\
(i)  $\mathscr G^{\partial}(E,D)$ contains  finitely many homotopy types, and  $\mathscr G(E,D)$ only   finitely many  diffeomorphism types; \\
(ii) $\mathscr C (E,D)$ contains  a finite  number of diffeomorphism types.
\end{corintro}

\noindent Notice that  we do not  bound a priori  the dimension of the manifolds in these classes:  this follows, in the aftermath, from bounding their entropy and diameter.

\noindent Moreover, notice that, while the results on Riemannian $2$-orbifolds and  ramified coverings (though concerning metrics of any possible sign) still pertain to the framework of spaces of negative curvature,  the classes   of non-geometric 3-manifolds and of high dimensional graph or cusp decomposable manifolds escape  from the realm of non-positive curvature. This is clear for non-prime 3-manifolds, and follows   from Leeb's work \cite{leeb} for irreducible 3-manifolds.  Cusp decomposable manifolds obtained by gluing boundary infra-nilmanifolds which are not tori do not admit nonpositively curved metrics, by the Solvable Subgroup Theorem; also, in  \cite{FLS} there are examples  of $n$-dimensional graph manifolds not supporting any locally $CAT(0)$-metric, for $n \geq 4$.


\vspace{3mm}

\small
{\em 
{\sc Acknowledgements}. The first author acknowledges financial support by  the Max-Planck Institut f\"ur Mathematik and praises the excellent working conditions provided by the Institut. Both the authors wish to thank R. Coulon for useful discussions during his stay in Rome in 2016, and S. Gallot for many precious hints.}
\normalsize

\vspace{3mm}
\section{Groups with acylindrical splittings}
\label{partI}

\small 
{\em 
\noindent This first part is  devoted to introducing some basic facts about groups with acylindrical splittings (Section \S\ref{subsecbasic}), and to developing  the combinatorial tools needed to prove the Entropy-Cardinality inequality (\ref{eqEC}) (Sections \S\ref{sspairwise}-\ref{secproofEC}).
}
\normalsize

\vspace{1mm}
  A finitely generated group $G$ possesses an {\em acylindrical splitting} if $G$ is  a  non-trivial amalgamated product  or HNN-extension, thus isomorphic to the fundamental group of a  non-trivial graph of groups $\mathscr G$,  and the canonical action of $G$ on the Bass-Serre tree associated to $\mathscr G$ is {\em acylindrical}. We recall the notion of \textit{acylindrical action} on a simplicial tree:

\begin{defn}[{\em Acylindrical actions on  trees}] \label{defacy}${}$\\
Let $G$ be a discrete group acting by automorphisms without edge inversions \linebreak (i.e. with no element  swapping the vertices of some edge) on a simplicial tree $\mathcal T$, endowed with its natural simplicial distance  $d_{\mathcal T}$ (i.e. with all edges of unit length).\linebreak
We shall say that the action 
$G\curvearrowright \mathcal T$ is \textit{$k$-acylindrical} if the set of fixed points of every elliptic element $g\in G$ has diameter less or equal to $k$, and that the action is \textit{acylindrical} if it is $k$-acylindrical for some $k$.
\end{defn}

 The notion of $k$-acylindrical action of a discrete group on a tree is a geometric reformulation of the notion of {\it $k$-step malnormal amalgamated product}, cp. \cite{karso}. 
 We recall that an element $g$ belonging to an amalgamated product $G=A\ast_C B$ is written in {\em normal form} when it is expressed  as $g= g_0g_1 g_2 \cdots g_n$, where $g_0 \in C$ and, for $i \geq 1$, no $g_i$ belongs to $C$   and two successive $g_i$ belong to different factors of the product; the integer $n=\ell(g)$ is then called the {\em syllable length} of $g$ (then, the  identity  and the elements of $C$ have  zero syllable length  by definition).\\
An amalgamated product $G=A \ast_C B$ is called \textit{$k$-step malnormal} if  $gCg^{-1}\cap C = \{1\}$  for all $g\in G$ with  $\ell(g) \geq k+1$ 
  (in particular, $k=0$  if and only if $C$ is a malnormal subgroup of $G$; and, by definition,  free products are $(-1)$-malnormal). \\
  A similar definition can be  given for a group $A\ast_\varphi =   \langle A, t\,|\,\mathrm{rel}(A), t^{-1}ct =\varphi(c) \rangle$ which is a HNN-extension of $A$ with respect to an isomorphism  $\varphi: C_- \rightarrow   C_+$ between subgroups $C_-, C_+$.  
Namely, by Britton's Lemma every element  $g\in G^\ast$ can be written in a {\em normal form} as 
   $$g=g_0t^{\varepsilon_1} g_1\cdots g_{m-1}t^{\varepsilon_m} g_m$$
where $g_0 \in A$, $\varepsilon_i =\pm1$,  and $g_i\in A\setminus C_{\varepsilon_i}$ if $\varepsilon_{i+1}=-\varepsilon_i$;
the syllable length of $g$ is defined in this case as $\ell(g)= m$.
Then, a  HNN-extension $G=A\ast_\varphi $ is called {\em $k$-step malnormal}  if   $g\,C_{\varepsilon}\,g^{-1}\cap C_{\varepsilon'} = \{1\}$ for  any $\varepsilon, \varepsilon' \in \{\pm\}$ and  for all $g\in G$ 
 with $\ell (g) \geq k+2$.
 \footnote{This seems to contrast with the definition for amalgamated products, but it actually yields that $G=A\ast_\varphi $ is a $0$-malnormal HNN extension  if and only if $C_+, C_-$ are  malnormal and conjugately separated in $A$,  and that $G$ is  $(-1)$-malnormal iff $C_\pm =\{1\}$. This is due to the relation $t^{-1} C_- t =C^+$; one might express the same condition by imposing  $g\,C_{\varepsilon}\,g^{-1}\cap C_{\varepsilon'} =\{1\}$ for all $g$ with $\ell (g) = k+1$ and whose normal form  satisfies some additional (awkward to write) restrictions.}

It is then easy  to check that a group $G$  which   is a  $(k-1)$-step malnormal amalgamated product or HNN-extension  admits a $k$-acylindrical action on his Bass-Serre tree. Conversely, if $G$ admits a $k$-acylindrical splitting   having a segment (resp. a loop) as underlying graph $\mathscr G$, then $G$ is a $(k-1)$-step malnormal amalgamated product (resp.  HNN-extension); see \cite{cer_pi1} for further details. 


\vspace{1mm}
  In what follows, we will be  interested in {\em non-elementary} acylindrical splittings: that, is, splittings of $G$ as a  non-trivial  amalgamated or HNN-extension  for which the action on the corresponding  Bass-Serre tree is non-elementary. 
Since non-trivial amalgamated products and HNN extensions act without global fixed points on their Bass-Serre tree,   this will be equivalent to asking that the action of $G$ is not {\em linear} (see the next subsection for  the basic terminology for acylindrical group actions on trees),  i.e. $G$  splits   as a {\em malnormal,  non-trivial}  amalgamated product $A\ast_C B$ (resp. HNN-extension $A \ast_\phi$)  {\em with $C$ of index greater than $2$ in  $A$ or $B$} \linebreak ({\em resp.} $[A:C_-]+[A:C_+ ]\ge 3$).

%
  

\vspace{5mm}
\subsection{Basic facts on acylindrical actions on trees} \label{subsecbasic}${}$ 

\noindent In the  following,  we will always   tacitly assume that the action of $G$  on the  simplicial tree $\mathcal T$ under consideration  is  by automorphisms and   without  edge inversions.\footnote{An action of a group $G$ on a simplicial tree $\mathcal T$ is said to be without (edge) inversions if for any oriented edge $\mathsf e$ and for any $g\in G$ we have $g\mathsf e\neq\bar{\mathsf e}$, where  $\bar{\mathsf e}$ is the  edge $\mathsf e$ with the opposite orientation. Notice that if a group $G$ acts on a simplicial tree $\mathcal T$ transitively on the (un)oriented edges, then $G$ acts without inversions on the barycentric subdivision of $\mathcal T$.} \linebreak
Then, the  elements of $G$ can be divided into two classes, according to their action: elliptic and hyperbolic elements. They are distinguished by their \textit{translation length}, which is defined,  for  $g\in G$,  as 
$$\tau(g)=\inf_{\mathsf v\in\mathcal T}d_{\mathcal T}(\mathsf v,g\cdot\mathsf v)$$
where $d_{\mathcal T}$ is the natural simplicial distance of ${\mathcal T}$, i.e. with all edges of unit length.
\linebreak
If $\tau(g)=0$ the element $g$ is called \textit{elliptic}, otherwise it is called \textit{hyperbolic}. \\
We shall denote by $\fix(g)$ the set of fixed points of an elliptic element $g$.  
We recall that $\fix(g)$ is a (possibly empty) connected subtrees of $\mathcal T$. 
If $h$ is a hyperbolic element then $\fix(h)=\varnothing$ and $h$  has a unique axis  
$$\ax(h)=\{\mathsf v\in\mathcal T\,|, d_{\mathcal T}(\mathsf v, h(\mathsf v))=\tau(h)\}$$
 on which it acts by translation: each element on the axis of $h$ is translated at distance $\tau(h)$ along  the axis, whereas elements at distance $\ell$ from the axis are translated of $\tau(h)+2\ell$. 

\noindent  Let $\mathcal T_G$  be the minimal subtree of  $\mathcal T$ which is $G$-invariant: the action of  $G$   is said {\em elliptic} it  $\mathcal T_G$  is a point, and {\em linear} if $\mathcal T_G$ a line; in both cases we shall say that the action of $G$ is {\em elementary}.

  \medskip

\noindent The next lemma resumes some facts about  centralizers and normalizers  of hyperbolic elements of acylindrical actions:

\begin{lem}\label{centr_norm_hyp}
	Let $G\curvearrowright \mathcal T$ be a $k$-acylindrical action on a  tree. \\
	Let $h\in G$ be hyperbolic, and let  $h_0$ be a hyperbolic element with  $\ax(h_0)=\ax(h)$ and $\tau(h_0)$ minimal:
\begin{itemize}
\item[(i)]  any hyperbolic element $g$ with  $\ax(g)=\ax(h)$  is a multiple of $h_0$; 
\item[(ii)]  if $s\ax(h) = \ax(h)$ then  either $shs^{-1}= h$ or $ shs^{-1} = h^{-1}$, 
depending on whether $s$ preserves the orientation of $\ax(h)$ or not;
\item[(iii)] the centralizer  $Z (h)$ is the infinite cyclic subgroup generated by $h_0$;
\item[(iv)] the normalizer $N(h)$  is either equal to $Z(h)$, or to an infinite dihedral group generated by $\{h_0, \sigma\}$, where $\sigma$ is an elliptic element swapping $\ax(h_0)$;
\item[(v)] if the action is   linear then $G$ is virtually cyclic.

\end{itemize}
\end{lem}

{\em Proof}.
Since $h_0$ has   minimum  translation length among hyperbolic elements with the same axis, any hyperbolic  $g$ with  $\ax(g)=\ax(h)$ acts on $\ax(h_0)$ as $h_0^n$, for some $n \in \mathbb{Z}$; hence  $g=h_0^n$ necessarily, by acylindricity, which shows (i).  \\
To see (ii), notice that  if  $\ax(h)=s\ax(h)=\ax(shs^{-1})$  then $shs^{-1}=h^{\pm 1}$ by (i), since $\tau(shs^{-1})=\tau(h)$.
Moreover, as $s$ stabilizes $\ax(h)$, if  $s$ is hyperbolic  then it has the same axis as $h$, preserves the orientation of the axis  and commutes with $h$; that is $shs^{-1}=h$.  
On the other hand, if $s$ is elliptic, by acylindricity it can only swap  $\ax(h)$ by a reflection through a vertex (because $G$ acts withouts edge inversions),  so $shs^{-1}$ acts on $\ax(h)$ as $h^{-1}$; this implies that   $shs^{-1}=h^{-1}$, again by acylindricity.\\
Assertion (i) implies that the subset  of hyperbolic elements of $G$ whose axes are equal to $\ax(h)$ is included in $Z(h)$.  Reciprocally, if $s\in Z(h)$ then, by assumption, $shs^{-1}=h$ and $s\ax(h)=\ax(h)$; therefore,  by the above discussion,
$s$ is a hyperbolic element which has the same axis as $h$. This proves (iii). \\
To see (iv), Let  $s\in N(h)$. As   $s\langle h_0\rangle s^{-1}=\langle h_0\rangle$ and conjugation is an automorphism of $G$, we deduce that $sh_0s^{-1}=h_0^{\pm 1}$, so $s$ stabilizes $\ax(h)$.
If $s$ is hyperbolic, we know by (ii) that it belongs to the infinite cyclic subgroup $\langle h_0\rangle$.  On the other hand, we have seen that  if $s$  is elliptic then it acts on $\ax(h)$ by a reflection through  a vertex $\mathsf v$, so  $s^2=1$ and $sh_0s^{-1}=h_0^{-1}$ by acylindricity. If $s'$ is another element of $N(h)\setminus Z(h)$, then also $s'$ acts by swapping $\ax(h)$ and fixing some $\mathsf v' \in \ax(h)$. Then, $s's^{-1}=s's$ acts on $\ax(h)$ as $h_0^n$, and by acylindricity $s's=h_0^{n}$. This shows that, if $N(h)\neq Z(h)$ then $N(h)\cong D_\infty$.

\noindent Finally, let us show (v). Assume that the minimal $G$-invariant  subtree of  $\mathcal T$ is a line. Since the action is non-elliptic, there exists   a  hyperbolic element $h$ preserving this line 
(as $G$ acts on  $\mathcal T$ without global fixed points, its elements cannot be all elliptic, or there would exist 
$s_1, s_2$ with  $\fix(s_1)\cap \fix(s_2)=\varnothing$ and this would produce a hyperbolic element $h=s_1s_2$).
Any other element $s \in G$ leaves $\ax(h)$ invariant, so  it is in $N(h)$ by (ii), and then $G$ is cyclic or dihedral by  (iii) and (iv).  
\qed\\

 The next lemma bounds  the size of the  intersection of the axes of two hyperbolic elements  in terms of the acylindricity constant $k$ and of the translation lengths: 

\begin{lem}\label{max_transl>diam}
	Let $G\curvearrowright \mathcal T$ be a non-elementary, $k$-acylindrical action on a  tree. \\
	 Let $h_1, h_2\in G$ be two hyperbolic elements with distinct axes: then,
	$$ \diam(\ax(h_1)\cap\ax(h_2)) \le \tau(h_1) + \tau(h_2) + k$$
\end{lem}

{\em Proof}. If $\ax(h_1)\cap\ax(h_2)=\varnothing$ the assertion is trivially verified.
Up to taking inverses, we may assume that the elements $h_1$, $h_2$ translate $I\!=\!\ax(h_1) \cap \ax(h_2)\neq\varnothing$ \linebreak in the same direction, and we let $I$ be oriented by this direction. 
We argue by contradiction: assume that $\diam(I) \ge \tau(h_1)+\tau(h_2)+k+1$.
Then, the element $[h_1,h_2]$ would fix pointwise the initial subsegment of length $k+1$ of $I$ and this, by $k$-acylindricity, implies that  $[h_1,h_2]=1$. This contradicts the assumption that  $h_1$ and $h_2$  have distinct axes.
\qed

\vspace{5mm}
\subsection{Schottky and pairwise-Schottky subgroups}
\label{sspairwise} ${}$ 

\noindent Recall that the subgroup $\langle g_1,...., g_n\rangle$ of $G$ generated by  $g_1,...,g_n$  is a \textit{Schottky subgroup} (of rank $n$), if for any $i=1,.., n$ it is possible to find subsets $X_{i}\subset\mathcal T$ for $i=1,..., n$ such that:
\begin{itemize}
	\item[(i)] $X_i\cap X_j=\varnothing$ for $i\neq j$;
	\item[(ii)]  $g_i^{\pm 1}(\mathcal T\smallsetminus X_i)\subseteq X_i$ for all $i$.
\end{itemize}
The $X_i$'s are called \textit{joint domains of attraction} of the $g_i$'s.

\begin{lem}\label{2schottky}
	Let $G$ be a group acting  on a simplicial tree $\mathcal T$. Let $h_1, h_2$ be two hyperbolic elements such that:
	$$\diam(\ax(h_1)\cap\ax(h_2))<\min\{\tau(h_1)\,;\,\tau(h_2)\}$$
	Then the group $\langle h_1, h_2\rangle$ is a rank $2$ free, Schottky  subgroup  of $G$.
\end{lem}

	This lemma is folklore (see for instance \cite{kaweid}). However let us describe  the domains of attraction $X_1$, $X_2$  for future reference. Let $p_i:\mathcal T\f\ax(h_i)$ be the projection  and  define $J_1=p_1(\ax(h_2))$ and $J_2=p_2(\ax(h_1))$: then one can choose $X_{i}=\{\mathsf v\in\mathcal T\,|\,p_i(\mathsf v)\not\in J_i\}$.

\begin{defn}[pS-family]\label{pS}
Let $G$ act on a simplicial tree $\mathcal T$  and let $\mathcal H=\{h_1,..., h_m\}$ be a collection of hyperbolic elements such that:
		$$\diam(\ax(h_i)\cap\ax(h_j))<\min\{\tau(h_i), \tau(h_j)\}, \quad\mbox{\small for any \normalsize}i\neq j$$
		We shall call $\mathcal H$ a {\it pairwise Schottky family} (shortly, pS-family) and let  $X_{i,j}$, $X_{j,i}$  be the domains of the pair $\{h_i, h_j\}$, $1\le i<j\le m$ (namely, where $X_{i,j}$ is the domain of attraction of $h_i$, with respect to the pair $\{h_i, h_j\}$).
\end{defn}

\begin{rmk}
	It is worth to stress the difference between a Schottky subgroup and the group $H$ generated by a pS-family $\mathcal H=\{h_1,..., h_n\}$. By definition, for any pair of elements $h_i, h_j\in\mathcal H$ of a pS-family it is possible to find   domains of attraction $X_{i,j}$ and $X_{j,i}$, so that $\{ h_i, h_j\}$ are in Schottky position and generate a rank $2$ free subgroup of $G$. However the whole collection $\mathcal H$ generally is not in  Schottky position, as the domains $X_{i,j}$, $X_{j,i}$ are not joint domains of attraction for all the $h_i$'s, as they depend on the pair $\{h_i, h_j\}$ (notice that the intersection $X_i:=\bigcap_{j\neq i} X_{i, j}$ might be empty for some index $i$). In particular $H$ generally is not free.
\end{rmk}

\vspace{3mm}
\subsection{Large pS-families with universally bounded $S$-lengths}
\label{sslarge} ${}$ 

\noindent We  start now considering finitely generated {\em marked groups} $(G,S)$ acting on trees, where  $S$ is any finite generating set for $G$. We will denote  by $| \; |_S$ the associated  word metric, and   by $S^n$ the relative ball of radius $n$, centered at the identity element. \linebreak
Let us also denote by $S_{hyp}^n$  the subset of hyperbolic  elements in $G$ of $S$-length smaller than or equal to $n$. 

 We start noticing  that, when the action is non-elementary, there always exist hyperbolic elements of $S$-length at most $2$ (see for instance \cite{CS}, \S2):
 
\begin{lem}
\label{gensys}
	Let $G\curvearrowright \mathcal T$ be an action on a simplicial tree, and let $S$ be any finite generating set for $G$:
\begin{itemize}

\item[(i)] if the action is non-elliptic, then there exists a hyperbolic element $h \in  G$ such that $|h|_S\leq 2$. Namely, either $h \in S$, or $h$ is the product of two elliptic elements $s_1,s_2 \in S$ with $\fix(s_1) \cap  \fix(s_2) = \emptyset$;	
				
\item[(ii)] if the action is non-elementary, then for any hyperbolic element $h \in G$ there exists $s \in S$ which does not belong to the normalizer $N(h)$.	
			
		
\end{itemize}
\end{lem}

 The next two Propositions are preparatory to the Entropy-Cardinality inequality.  \linebreak
They  show  that, from any generating set $S$ of a group $G$ acting acylindrically and non-elementarily on a tree, one can always produce a large (compared to $|S|$) pS-family $\mathcal H$ of hyperbolic elements with universally bounded $S$-length,   whose axes are, moreover, in one of two basic configurations.


\begin{prop}[Large conjugacy classes of  hyperbolics  with bounded $S$-lengths]\label{tanti_ip_assi_dist}
	Let $G\curvearrowright \mathcal T$ be a non-elementary, $k$-acylindrical action on a tree, and let $S$ be any finite generating set for $G$.
	There exists a subset ${\mathcal H} \subset S_{hyp}^{20}$ of hyperbolic elements, with $|{\mathcal H}|  \geq \sqrt[4]{|S|}$ such that:
	
\begin{itemize}
\item[(i)]  the elements of ${\mathcal H}$ have all different axes;
\item[(ii)]  every element of ${\mathcal H}$ has at least $\sqrt[8]{|S|} /2$ conjugates with distinct axes in $S_{hyp}^{20}$.
\end{itemize}	
\end{prop}	

  \pagebreak
	
\begin{prop}[Large pS-families with universally bounded $S$-lengths]\label{unif_bounded_large_pS-family}\hfill\\	
      Let $G\curvearrowright \mathcal T$ be a non-elementary, $k$-acylindrical action on a tree.  For any finite generating set $S$ of $G$ there exists a pairwise Schottky family $\mathcal S$   with the following properties:
	\begin{itemize}
		\item[(i)] all the elements of $\mathcal S$ have the same translation length;
		\item[(ii)] $|h|_S\le 20(k+3)$ for all $h \in \mathcal S$;
		\item[(iii)] $|\mathcal S|\ge \sqrt[16]{|S|} /\sqrt 2$.
		\item[(iv)] either there exists $h_1\in\mathcal S$ such that $\ax(h)\cap\ax(h_1)\neq\varnothing$ for all $h\in\mathcal S$, \\or all the axes of the elements in $\mathcal S$ have pairwise empty intersection.
	\end{itemize}
\end{prop}

In order to prove the above Propositions,  we preliminary  show:

\begin{lem}\label{tanti_ip_distinti}
	 Under the   above assumptions, we have 
	$|S_{hyp}^4|\ge\sqrt{|S|}$. Moreover, in $S_{hyp}^4$ there are at least two hyperbolic elements having distinct axes.
	\end{lem}

 {\em Proof of Lemma \ref{tanti_ip_distinti}.}
 Let $h\in G$ be a hyperbolic element such that $|h|_S\le 2$.
 Consider the set $C_S(h) \subset S^4 $ of conjugates $shs^{-1}$ of $h$ by all the elements of $S$. 
 Observe that either $|C_S(h)|\ge\sqrt{|S|}$, or there exists a subset $S_h \subset S$ of cardinality at least $\sqrt{|S|}$ such that $shs^{-1}=s'hs'^{-1}$ for all $s,s' \in S_h$. In the first case we deduce that  $|S_{hyp}^4| \ge\sqrt{|S|}$. In  the second case, choose $s_0 \in S_h$ and consider the map 
 $F : S_h \f S^2$ defined by $ F(s)=s^{-1}s_0$.	
 By definition, $F$ is an injective map.  Moreover, since $s_0hs_0^{-1}=shs^{-1}$ for all $s \in S_h$, we have 
 $s^{-1}s_0 h =hs^{-1}s_0$, so $\mathrm{Im}(F)\subset Z(h)$. By Lemma \ref{centr_norm_hyp} we know that $Z(h)$ is infinite cyclic, generated by a hyperbolic element $h_0$;   we deduce that $\mathrm{Im}(F)\subset S_{hyp}^2$, so $|S_{hyp}^4|\ge|S_{hyp}^2|\ge\sqrt{|S|}$.
	Finally, notice that, in any case,  $C_S(h)$ contains at least two hyperbolic elements with distinct axes: otherwise, every $s\in S$ would be in the  normalizer $N( h)$ and the action of $G$ on $\mathcal T$ would be elementary.\qed

\vspace{5mm}
  {\em Proof of Proposition \ref{tanti_ip_assi_dist}.}
 Consider the equivalence relation between hyperbolic elements $h,h'\in G$ defined by 
  $h\sim h'$ if and only if $\ax(h)=\ax(h')$, 
  and let $C_1,..., C_m$ be the equivalence classes in $S^4_{hyp}$ with respect to $\sim$.
 
\noindent By Lemma \ref{tanti_ip_distinti} we know that  $|S_{hyp}^4|\ge\sqrt{|S|}$, so:

(a) either there exists an equivalence class, say $C_1$, such that $|C_1|\ge\sqrt[4]{|S|}$;  

(b) or $m\ge\sqrt[4]{|S|}$. 

 \noindent In the first case, notice that we have  $m \geq 2$, always by Lemma \ref{tanti_ip_distinti}. \linebreak 
 So, in case (a), let $C_1 = \{ h_1,...,h_n \}$ with  $n \ge\sqrt[4]{|S|}$, and  let $h \in C_2$: then, the set  ${\mathcal H}=\{h_i h h_i^{-1} \; | \; h_i \in C_1 \}$ is a collection of $n$ hyperbolic elements in $S^{12}_{hyp}$ with distinct axes.
 Indeed, if $\ax(h_ihh_i^{-1})=\ax(h_jhh_j^{-1})$ then $h_{j}^{-1}h_i$ would preserve both $\ax (h_i)=\ax(h_j)$ and $\ax(h)$; this would imply $\ax(h)= \ax(h_i) = \ax (h_j)$, a contradiction.
Moreover, in this case there are also at least $n$ distinct conjugates (with distinct axes) of
 $h$ and of each  $h_i h h_i^{-1}$ in $S^{12}_{hyp}$, since the collection of elements $h_jhh_j^{-1}=h_jh_i^{-1}(h_i h h_i^{-1})h_ih_j^{-1}$ is also a collection of conjugates of $h_ihh_i^{-1}$.
Therefore,  in this case we conclude that 
$|\mathcal H|\ge n\ge\sqrt[4]{|S|}.$
 
\noindent In case (b), let $\mathcal R= \{ h_1,...,h_m \}$ be a set of representatives for the classes  $C_i$: we will now prove that every $h_i$ has at least $\sqrt{m}$ distinct conjugates 
in $S^{20}_{hyp}$.  \linebreak To this purpose,  let 
 $$ \mathcal C_{{\mathcal R}}(h_i) =\{  h_jh_i h_j^{-1} \; | \; h_j \in \mathcal R \}  \subset S^{12}_{hyp}$$ 
 Again, either $\mathcal C_{{\mathcal R}}(h_i)$ contains  at least $\sqrt{m}$ distinct conjugates of $h_i$ and we are done, or there exists a subset ${\mathcal H}_i \subset {\mathcal R}$ with cardinality $|{\mathcal H}_i | \geq \sqrt{m}$ such that
\begin{equation}
\label{hh'}
hh_i h^{-1}=h' h_i h'^{-1} =: g_i \mbox{ for all } h,h' \in {\mathcal H}_i
\end{equation}
In this case, consider the map  
$c_{h_i}: {\mathcal H}_i \rightarrow S^{20}_{hyp}$ 
defined by $h \mapsto c_{h_i}(h) := h^2 h_i h^{-2}$. To show that there are at least $\sqrt m$ distinct conjugates of $h_i$, it is enough to show that this map is injective. 
Actually, assume that $c_{h_i}(h)=c_{h_i}(h')$, for some $h,h' \in {\mathcal H}_i$: as, by (\ref{hh'}) we have that $h^{-1}h'=:\gamma_i \in Z(h_i)$, we would deduce
 $$g_i=hh_ih^{-1}= h^{-1}(h'^2h_ih'^{-2})h 
   = \gamma_i g_i \gamma_i^{-1}$$
that is $[g_i, \gamma_i]=0$.  Therefore,  $\gamma_i\in Z(h_i)\cap Z(g_i)=Z(h_i)\cap Z(hh_ih^{-1})$ which is possible if and only if $\ax(h_i)=h\,\ax(h_i)$. But then $h$ and $h_i$ would share the same axis, which contradicts the fact that  they are distinct elements of ${\mathcal R}$.\\
We have therefore showed that, in case (b),  every $h_i \in {\mathcal R}$  has at least  $\sqrt{m}$ distinct conjugates in  $S^{20}$.
Now,  if we have $\ax (sh_is^{-1}) = \ax (s'h_is'^{-1}) $ for $s,s' \in S$, then \linebreak 
$(s^{-1}s') \cdot \ax(h_i) = \ax (h_i)$,  so $(s^{-1}s')  h_i (s^{-1}s')^{-1}= h_i^\epsilon$  for some $\epsilon \in \{\pm 1\}$, by   Lemma \ref{centr_norm_hyp}; that is, either $s'h_is'^{-1}=sh_is^{-1}$, or $s'h_is'^{-1}= sh_i^{-1}s^{-1}  =  (sh_is^{-1})^{-1}$.
It readily follows that $h_i$ has at least $\sqrt{m} /2$ conjugates in $S^{20} $.\qed


\vspace{5mm}
  {\em Proof of Proposition \ref{unif_bounded_large_pS-family}.}
	By Proposition \ref{tanti_ip_assi_dist} we know that there is a collection $\mathcal H \subset S_{hyp}^{20}$ of hyperbolic elements which are all conjugates of one $h_0\in\mathcal H$, have distinct axes, and with cardinality  $ |\mathcal H |=n \ge \sqrt[8]{|S|}$. Moreover, for any pair of elements $h, h'\in\mathcal H $ we know  by Lemma \ref{max_transl>diam} that $$\diam(\ax(h)\cap\ax(h'))\le\tau(h)+\tau(h')+k=2\cdot\tau(h_0)+k$$  It follows that, for all $p>2+\frac{k}{\tau(h_0)}$:
	$$\diam(\ax(h^{p})\cap\ax(h'^{p}))<
	\tau(h^{p})=\tau(h'^{p})$$
	 for all $h, h'\in\mathcal H $. Therefore, the collection $\mathcal H_{1}=\{h^{k+3} \; | \; h \in \mathcal H\}$ provides a pairwise Schottky family satisfying (i), (ii) and (iii). We shall now prove that we can extract from $\mathcal H_1$ a subset $\mathcal S$ with $|\mathcal S|\ge\sqrt{n}$  satisfying one of the two conditions in (iv). 
Actually, let $m (h) $ be the number of elements of $\mathcal H_1$ whose axis  intersects $\ax(h)$, and let $m = \max_{h \in \mathcal H} m (h)$. 
Then either there exists  an element of $\mathcal H_1$, say $h_1$, with $m (h_1) \geq \sqrt{n}$, or   $m (h) < \sqrt{n}$ for all $h \in \mathcal H_1$.
In the first case, we define  $\mathcal S = \{ h \in \mathcal H_1 \;  | \; \ax(h) \cap \ax(h_1) \neq \emptyset \}$ and we are done.
Otherwise,  let $\mathcal S=\{ h_{1},..., h_{m} \} \subset \mathcal H_1$ be a maximal subset  of elements with pairwise disjoint axes, and let $\mathcal S_i \subset \mathcal H_1$ be the subset of elements whose axes intersect $\ax(h_i)$. \linebreak 
As $\mathcal S$ is maximal, we have
$$ \mathcal H_1 = \bigcup_{i=1}^m \mathcal S_i$$
so $n \leq m \cdot | \mathcal S_i | < m \sqrt{n}$. Hence, $| \mathcal S | = m > \sqrt{n}$.\qed
\vspace{2mm}


\vspace{4mm}
\subsection{Extraction of Schottky subgroups of large rank} \label{ssextraction} ${}$ 

\noindent In this section we consider a pS-family $\mathcal S=\{h_1,..., h_m\}$ with $\tau(h_i)=\tau$, $\forall i=1,...,m$ and whose axes satisfy one of the two configurations explained in Proposition \ref{unif_bounded_large_pS-family}. We will explain how to extract from the group $H=\langle\mathcal S\rangle$ generated by $\mathcal S$ a genuine Schottky subgroup of rank  $r \geq \frac{\sqrt{m}-1}{2}$, whose generators have $\mathcal S$-length   $\leq4$.



\smallskip
\subsubsection{pS-families distributed along one axis} 

\noindent  Here we study the case where all axes of the elements in $\mathcal S$ intersect an axis, say $\ax(h_1)$. We shall in particular examine two opposite (but not exhaustive) subcases:

\begin{itemize}
	\item[(a)] the case where {\em the intersection between the axes with $\ax(h_1)$ agglomerate together}, i.e. there exists a segment $I \subset\ax(h_1)$ with $\diam(I)\le 2\,\tau$ containing all the intersections $J_{1,i}:=\ax(h_1)\cap \ax(h_i)$ ;

	\item[(b)] the case where {\em the intersections with $\ax(h_1)$ are sparse along $\ax(h_1)$, } i.e.  $J_{1,i}=\ax(h_1)\cap\ax(h_i)\neq\varnothing$, $\forall i=2,...,M$ but $J_{i,j}:=\ax(h_i)\cap\ax(h_j)=\varnothing$ for every  $i, j\ge 2$ with $i\neq j$.
	
\end{itemize}
(Notice that $J_{1,i} = p_1 (\ax(h_i))$, so the notation is consistent with the one in  \S\ref{sspairwise}).
For easy future reference we will call such families {\em agglomerated} (in case (a)) or {\em sparse} (in case (b)) with respect to $h_1$.
\vspace{1mm}

 The next Proposition shows that these configurations are actually more frequent than one can expect:

	\begin{figure}[H]
		\centering
		\includegraphics[scale=0.30]{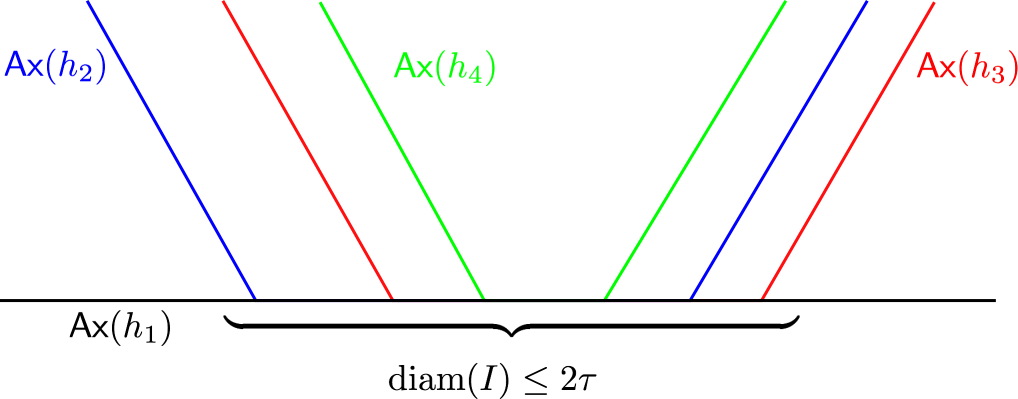}
		\caption{\small pS-family agglomerated with respect to $h_1$}
		\label{lemma36}
         \end{figure}

	\begin{figure}[H]
		\centering
		\includegraphics[scale=0.30]{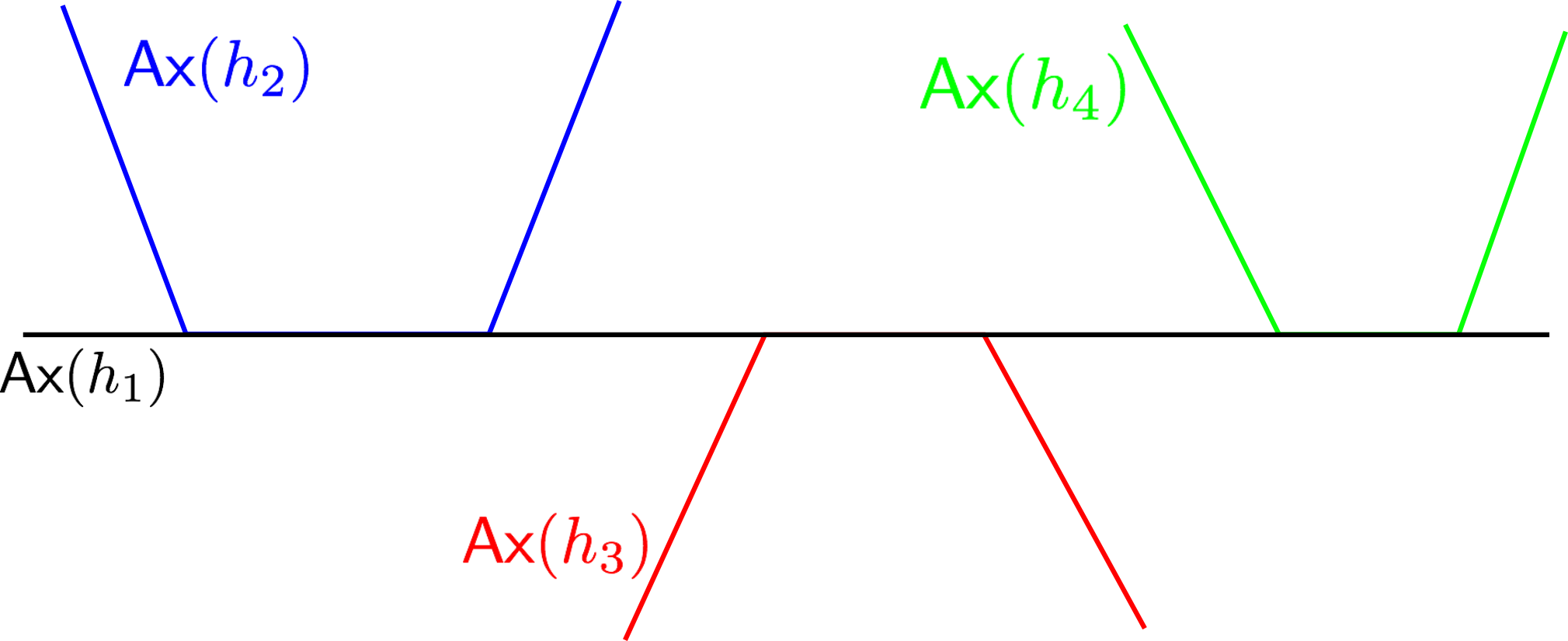}
		\caption{\small pS-family sparse with respect to $h_1$}
		\label{lemma36}
	\end{figure}

\begin{prop}\label{extraction_intersecting}
	Let $\mathcal S=\{h_1,..., h_m\}$
	be a pS-family with  $\ax(h_1)\cap\ax(h_i)\neq \varnothing$ for all $i$.
	Then, there exists a subset $\mathcal S' \subset \mathcal S$,  containing $h_1$,  which either is   agglomerated  with respect to $h_1$ and has  cardinality  $m' \ge \sqrt m$, or is sparse  with respect to $h_1$ and has cardinality $m' \ge \sqrt m +1$.
\end{prop}

{\em Proof}.
	For any $k=2,..., m$ let us denote as before  $J_{1,k}=\ax(h_1)\cap \ax(h_k)$ and let $\mathcal S^*=\mathcal S\smallsetminus\{h_1\}$. Since $\mathcal S$ is a pS-family such that $\tau(h_i)=\tau$, we know that  $\diam(J_{1,k})<\tau$. Orient $\ax(h_1)$ by the translation direction of $h_1$. Then, to each $J_{1,k}$ we associate $\mathsf v_{k}^-$, $\mathsf v_{k}^+$, respectively the initial and final vertex of $J_{1,k}$ with respect to the orientation of $\ax(h_1)$. We then consider the following order relation on $\{h_k\}_{k=2}^{m}$: we shall say that $h_k\le h_j$ if
	\begin{itemize}
		\item either $\mathsf v_k^-\le\mathsf v_j^-$;
		\item or $\mathsf v_k^+\le \mathsf v_j^+$ and  $\mathsf v_k^-=\mathsf v_j^-$;
		\item or $k\le j$ and $J_{1,k}=J_{1,j}$.
	\end{itemize}
	The relation $\le$ defines a total order on $\mathcal S^*$. Up to re-indexing the elements  in $\mathcal S^*$ we may assume that  $h_i<h_j$ if and only if $i<j$. Set $i_1=2$ and let $I_1$ be the interval on $\ax(h_1)$ starting at $\mathsf v_{i_1}^{-}$ and having length $2\tau$, and define $C_{I_1}$ be the set of those elements in $\mathcal S^*$ such that $J_{1,k}\subseteq I_1$ and $m_1=|C_{I_1}|$. By definition $C_{I_1}$ is non-empty since $h_{2}\in C_{I_1}$. If $m_1\ge\sqrt{m}-1$, then we define  $\mathcal S'=\{h_1\}\cup C_{I_1}$ and we are done, as $\mathcal S'$ is an agglomerated family with respect to $h_1$. Otherwise set $i_2=m_1+2$, consider the segment $I_{2}$ starting at $\mathsf v_{i_2}^{-}$ of length $2\tau$ and define $C_{I_2}=\{h_k\in \mathcal S^*\smallsetminus C_{I_1}\,|\, J_{1,k}\subseteq I_2\}$ and $m_2=|C_{I_2}|$. Then either $m_2\ge\sqrt{m}-1$, and we conclude taking $\mathcal S'=\{h_1\}\cup C_{I_2}$, or we set  $i_3=m_1+m_2+2$.  We proceed in this way until either we find a $C_{I_k}$ whose cardinality is greater or equal to $\sqrt{m}-1$, or we exhaust the set $\mathcal S^*$. Assume that the process ends after $K$ steps and that all of the $C_{I_k}$ have cardinality strictly smaller than $\sqrt{m}-1$: then, we consider the set $\{h_{i_1},...,h_{i_K}\}$.  Notice that, since $\diam(I_k)=2\tau$ and $\diam(J_{1,i})<\tau$ for all $i=2,...,m$ we deduce  that $d_{\mathcal T}(\ax(h_{i_k}),\ax(h_{i_j}))>0$ if $k\neq j$. Therefore, the collection $\{\ax(h_{i_1}),...,\ax(h_{i_K})\}$ is a collection of pairwise disjoint axes, such that each element in the collection has a non-empty intersection with $\ax(h_1)$; so $\mathcal S'=\{h_1\}\cup\{h_{i_1},...,h_{i_K}\}$ is a pS-family sparse with respect to $h_1$ and such that $\tau(h')=\tau$ for any $h'\in\mathcal S'$. Observe that
 $ K>\frac{m-1}{\sqrt m -1} \ge \sqrt m $, so $|\mathcal S'|\ge\sqrt m +1$.
\qed


\vspace{3mm}
We now show how to extract a free subgroup of rank  $r \geq m'-1$ from a pS-family $\mathcal S'=\{h_1,...,h_{m'}\}$ whose axes are in configuration (a) or (b).

\begin{lem}\label{pS_bigcap}
	Let $\mathcal S' =\{h_1,...,h_{m'} \}\subset G$ be a pS-family with $\tau(h_i)=\tau$ for all $i$,  agglomerated with respect to $h_1$. Then,  $\langle h_1^4,...,h_{m'}^4\rangle$ is a Schottky subgroup of $G$.
\end{lem}

{\em Proof}.  Let $I \subset\ax(h_1)$ be a segment with $\diam(I)\le 2\,\tau$ containing all the intersections $J_{1,i}=\ax(h_1)\cap \ax(h_i)$, let $J_{i,j}=\ax(h_i)\cap\ax(h_j)$ and let 
	 us denote by  $\mathcal T' \subset \mathcal T$  the subset: 
	\small
	$$\mathcal T' = I \cup\left(\bigcup_{i, j\ge 2,\,i\neq j} J_{i,j}\right)$$
	\normalsize
Since $J_{1,i}\neq\varnothing$ and is contained into $I$ for all $i \neq 1$, we conclude that $\mathcal T'$ is a connected subtree of $\mathcal T$. Moreover, since $\tau(h_i)=\tau$ for all $i$, by definition of agglomerated, pS-family we have
$\diam(\mathcal T')\le \diam(I)+2\,\max_{i\neq j}\diam(J_{i,j}) < 4 \tau$.\linebreak
It follows that the sets $X_i=\{\mathsf v\,|\, p_i(\mathsf v)\not\in\mathcal T' \}$ are joint domains of attraction for the collection $\{h_i^4\}_{i=1,...,m'}$. Actually, $X_i\cap X_j=\varnothing$ for any $i\neq j$  and if  $\mathsf v\in X\smallsetminus X_i$  the projection of $\mathsf v$ on $\ax(h_i)$ is contained in $\mathcal T'$, as $\tau(h_i^{\pm 4})=4\tau>\diam(\mathcal T')$; hence the projection of $h_i^{\pm 4}(\mathsf v)$ on $\ax(h_i)$ lies in $\ax(h_i)\smallsetminus \mathcal T'$, \textit{i.e.} $h_i^{\pm 4}(\mathsf v)\in X_i$.
\qed

\vspace{3mm}
\begin{lem}\label{pS-almostdisjoint}
	Let $\mathcal S'=\{h_1,..., h_{m'}\}\subset G$ be a pS-family with $\tau(h_i)=\tau$ for all $i$, which is sparse with respect to $h_1$. Then, $\langle h_2,..., h_{m'}\rangle$ is a Schottky subgroup of $G$.
\end{lem}

{\em Proof}.
For each $i=2,..,m'$   define $X_i=\{\mathsf v\in\mathcal T\,|\, p_i(\mathsf v)\not\in J_{1,i}\}$. By definition  the projections of the $X_i$'s on $\ax(h_1)$ are disjoint;  so,  if $\mathsf v\in X\smallsetminus X_i$, then   its projection on $\ax(h_i)$ lies in $J_{1,i}$. As $\tau(h_i)>\diam(J_{1,i})$ it follows that $h_i^{\pm1}(\mathsf v)\in X_i$.
\qed

\newpage

\vspace{-2mm}
\begin{figure}[h]
	\centering
	\includegraphics[scale = 0.21]{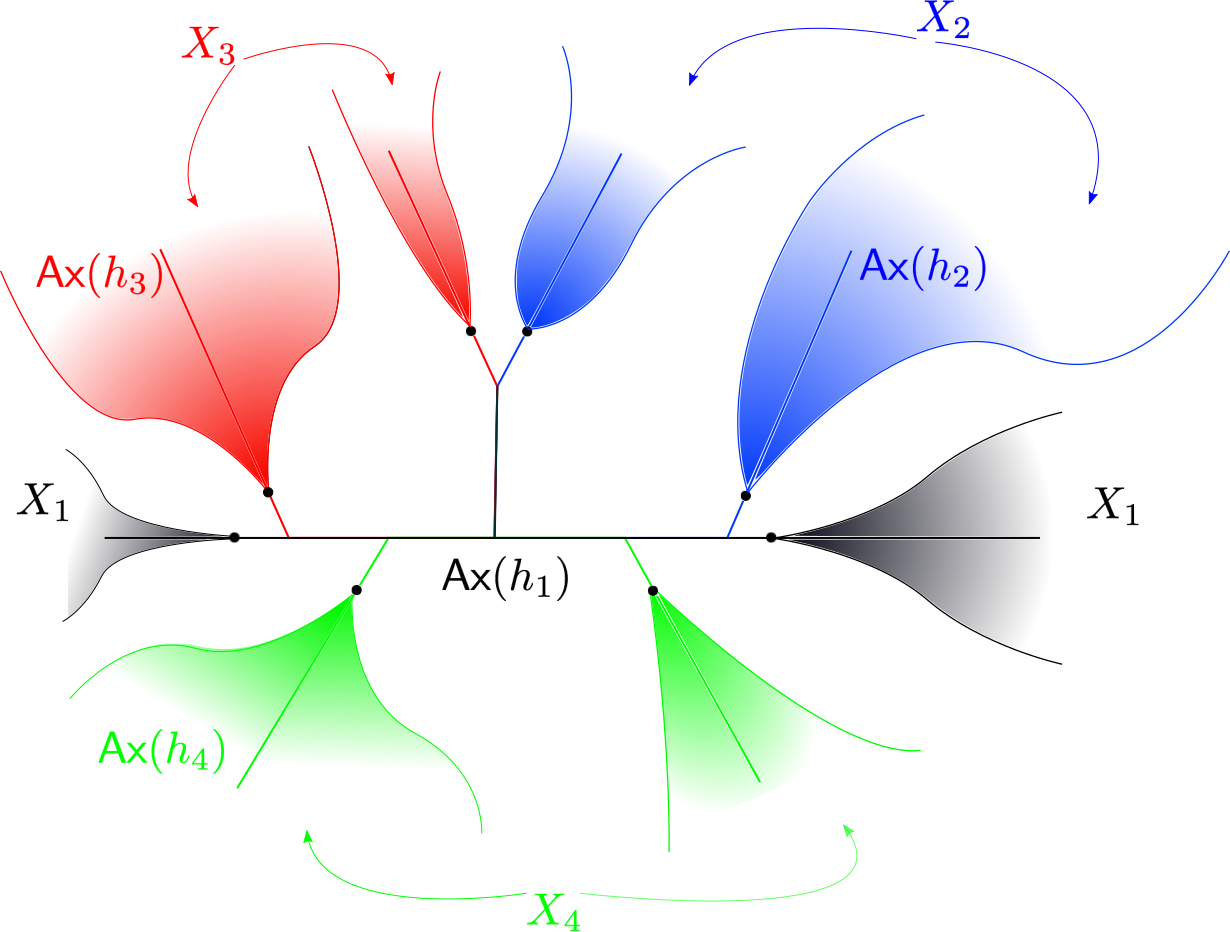}
	\caption{\vspace{20mm}\small Joint domains of attraction as in Lemma \ref{pS_bigcap}}
	\label{42}
\end{figure}
\vspace{-21mm}

 \subsubsection{pS-families with disjoint axes}  ${}$

\begin{figure}[H]
\includegraphics[scale = 0.27]{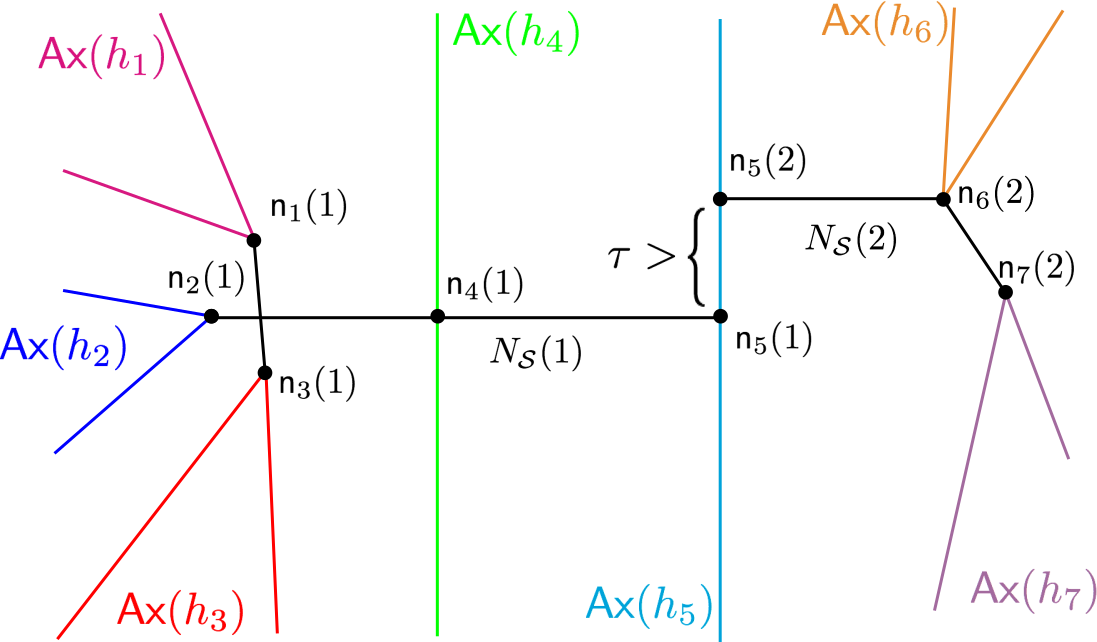}
\caption{\small Disjoint pS-family with small projections}
\label{nerves_nervertices}
\end{figure}

\begin{figure}[H]
	\centering
\includegraphics[scale = 0.27]{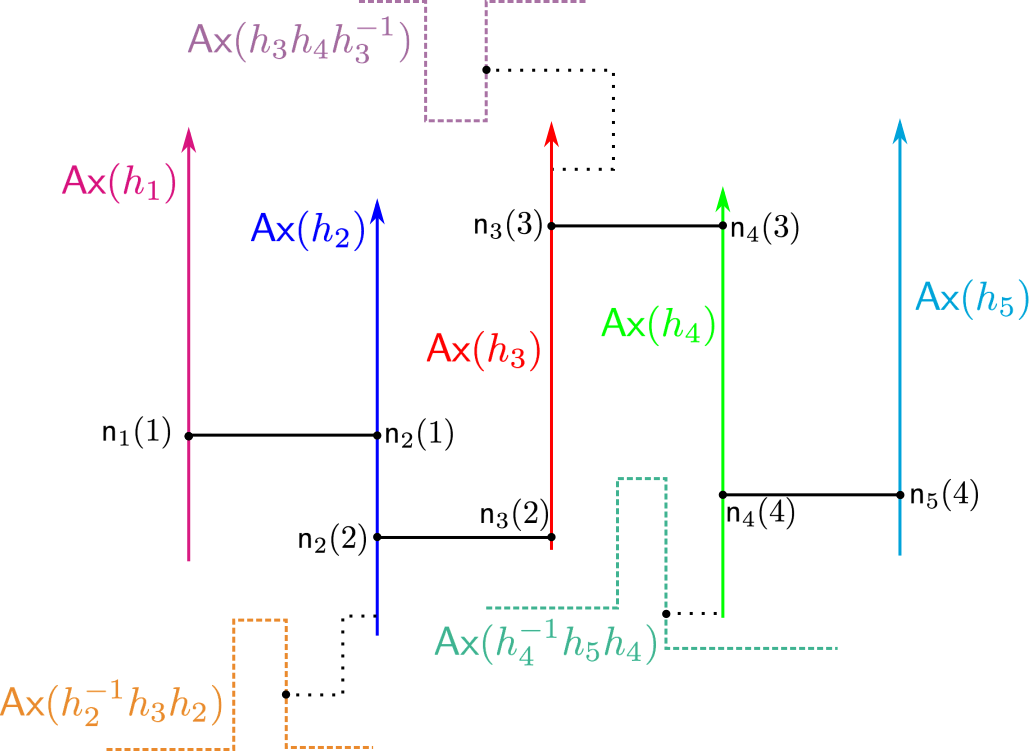}
\caption{\small Disjoint sequential pS-family (and a new one with one nerve)}
\label{schottky_from_disjoint_paths}
\end{figure}

\noindent  We  now study  the case where the axes  of the  family $\mathcal S=\{h_1,..., h_m\}$ are all disjoint.
 To properly describe this case, let us introduce  some notation. Let $p_i:\mathcal T\f \ax(h_i)$ \linebreak
 be the projection, let $\mathcal A=\bigcup_{i=1}^m\ax(h_i)$ and let $\mathcal T_{\mathcal S}$ be the convex hull of $\mathcal A$. \linebreak
 We shall call the \textit{nerves} of $\mathcal S$ the connected components $N_{\mathcal S}(j)$ of $\overline{\mathcal T_{\mathcal S}\setminus\mathcal A}$.   \linebreak
 We remark that, if $\ax(h_i)\cap N_{\mathcal S}(j)\neq \varnothing$, then the intersection consists of a single vertex  $\mathsf n_{i}(j)$, which coincides with the projection $p_i(\ax(h_k))$ for each $\ax(h_k)$ whose intersection with $N_{\mathcal S}(j)$ is non-empty; we  shall call it the \textit{i-nervertex of $N_{\mathcal S}(j)$}.\\
We will focus on 2 particular subcases:
\begin{itemize}
	\item[(c)] the case where for each  $ i=1,...,m$ the subset  $J_i= \bigcup_{j \neq i} \{p_i(\ax(h_j)) \; |\; j \neq i \}$ of the projections on $\ax(h_i)$ of all the other axes is included in a subsegment $I_i\subset\ax(h_i)$ of length strictly less than $\tau$  (see  Fig.\ref{nerves_nervertices});
	\item[(d)] the case where there are precisely $m-1$ nerves $N_{\mathcal A}(1)$,..., $N_{\mathcal A}(m-1)$ and for any $j=1,.., m-1$ the nerve $N_{\mathcal A}(j)$ coincides with the geodesic segment connecting $\ax(h_j)$ to $\ax(h_{j+1})$ (see Figure \ref{schottky_from_disjoint_paths}, disregarding the dotted lines and the  $\ax(h_ih_4h_i^{-1})$'s).
\end{itemize}
We shall refer to such configurations as to a {\em disjoint family with small projections}
(in case (c)) and a {\em disjoint sequential family} (in case (d)). 
   \vspace{-2mm}

\pagebreak

   \vspace{-5mm}
\noindent Again, the next Proposition shows that  the configurations (c) and (d) are frequent.

 \begin{prop}\label{disjoint_schottky}
 	Let $\mathcal S=\{h_1,...,h_m\}$  be a  pS-family with pairwise disjoint axes. Then,  there exists  a subset $\mathcal S'\subset\mathcal S$  
	of cardinality $m'\ge\sqrt{m}$ which is either a disjoint pS-family with  small projections   or  a disjoint sequential  pS-family.
 \end{prop}

\noindent To better handle  in the proof of Proposition \ref{disjoint_schottky} the possible configurations of the axes of a general pS-family $\mathcal S=\{h_1,...,h_m\}$ with pairwise disjoint axes, let us  define the  {\em graph $\Gamma_{\mathcal S}$  of the nerves of ${\mathcal S}$}, as follows: its vertex set $\mathsf V(\Gamma_{\mathcal S})$ is given by the nerves $N_{\mathcal S}(j)$ and by the axes $\ax(h_i)$ of the elements of $\mathcal S$, and there is an edge between $N_{\mathcal S}(j)$ and $\ax(h_i)$ if and only if $N_{\mathcal S}(j)$ has a $i$-nervertex $\mathsf n_{i}(j)$, i.e. $N_{\mathcal S}(j)\cap\ax(h_i)\neq\varnothing$. 	
   \vspace{-2mm}
 
 \begin{figure}[H]
 	\centering
 	\includegraphics[scale=0.25]{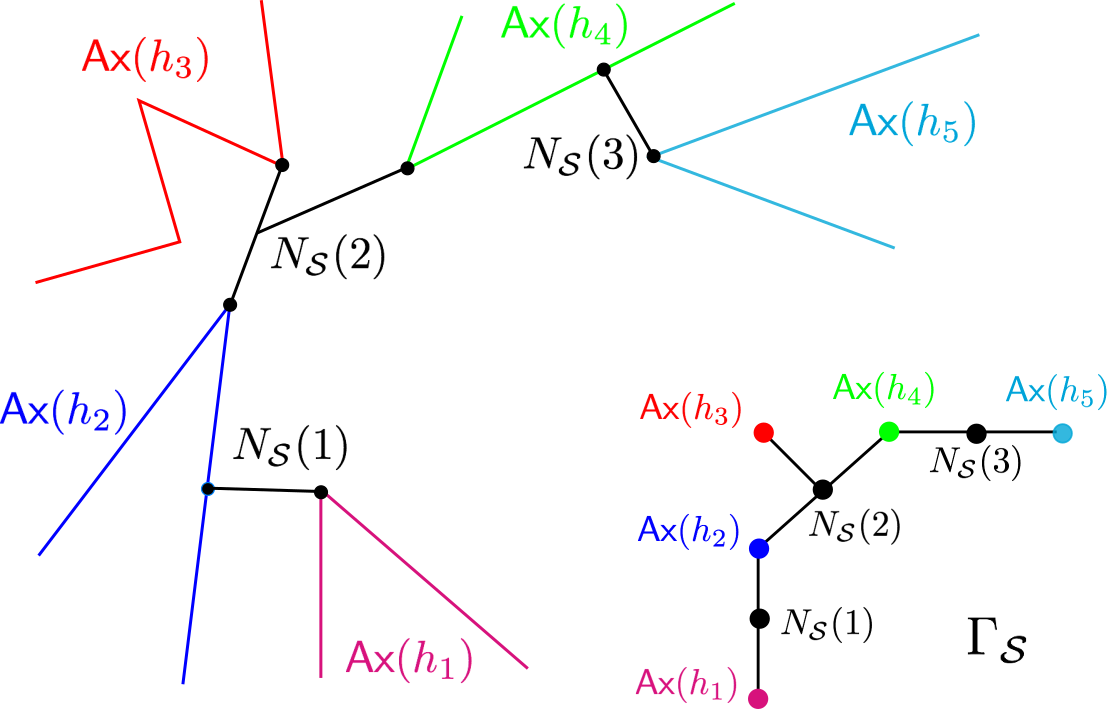}
 	\caption{The graph of nerves}
 	\label{networkgraph}
  \end{figure}

   \vspace{-5mm}
 \noindent Notice that, as it is build from the tree $\mathcal T_{\mathcal S}$, the graph $\Gamma_{\mathcal S}$  is a tree. The leaves of $\Gamma_{\mathcal S}$ are all vertices of type $\ax(h_i)$, for suitable $i$'s, and  there are no edges connecting  two vertices of type $N_{\mathcal S}(j)$ or two vertices of type $\ax(h_i)$. Moreover, we clearly have  
 \small
 \begin{equation}\label{eqvertices}
 m+1\le |\mathsf V(\Gamma_{\mathcal S})|\le 2m-1
\end{equation}
 \normalsize
with $|\mathsf V(\Gamma_{\mathcal S})|=m+1$ implying that $\mathcal S$ is a disjoint family with small projections (namely, a single point $n_i(j)$ on each axis $\ax(h_i)$) and $|\mathsf V(\Gamma_{\mathcal S})|=2m-1$ if and only if the degree of any vertex of type $N_{\mathcal S}(j)$ is equal to $2$.  

 \begin{lem} \label{tree}
 Let $\Gamma$ be a connected tree with $n$ vertices, and let $\ell$, $d$ denote respectively the number of leaves and the diameter: 
 if $d$ is even,    then $\ell \cdot d \geq 2(n-1)$.
  \end{lem}

  {\em Proof of Lemma  \ref{tree}.}
  Let $\gamma$ be a diameter of $\Gamma$, with central point $\mathsf p$ (which is  a vertex, since $d$ is even). 
  Summing the number of  vertices different from $\mathsf p$ belonging to  all segments $\overline{\mathsf p\mathsf v_i}$, where $\mathsf v_i$ runs over the leaves of  $\Gamma$, we find
  $\ell \cdot \frac{d}{2}  \geq n-1$, which is the announced inequality.\qed
   \vspace{4mm}

  {\em Proof of Proposition \ref{disjoint_schottky}.}
 	Let $n$ be the number of vertices of the tree  $\Gamma_{\mathcal S}$.
Notice that, by the very  definition of $\Gamma_{\mathcal S}$,    the diameter $d$ is realized by a geodesic $c$  connecting two leaves and is even.
Hence, by the above lemma,  either the number of leaves of    $\Gamma_{\mathcal S}$ is $\ell\ge\sqrt{(n-1)}$  or its diameter is  $d\ge2\sqrt{(n-1)}$. Assume that the first case holds:  each leaf of $\Gamma_{\mathcal S}$ corresponds to an axis $\ax(h_i)$ so let $\mathcal S'\subset \mathcal S$ be the subset of those elements whose axes correspond to the leaves of $\Gamma_{\mathcal S}$. For any $h, h_i\in\mathcal S'$ the projection of $\ax(h)$ on $\ax(h_i)$ is one point, namely the nervertex $\mathsf n_{i}(j)$ of the only nerve $N_{\mathcal S}(j)$ connected to $\ax(h_i)$. Thus,  the collection $\mathcal S'$ is a disjoint pS-subfamily of $\mathcal S$ with small projections and 
  $|\mathcal S'|=\sqrt{(n-1)}\ge\sqrt{m}$, by (\ref{eqvertices}). 
Assume now that we are in the second case. Let $\ax(h_{j_1}),\dots, \ax(h_{j_{d+1}})$ be the axes corresponding to the odd vertices of the geodesic $c$: then, the set $\mathcal S'=\{h_{j_1},... h_{j_{p+1}}\}$ is  a sequential disjoint pS-subfamily of $\mathcal S$ with 
$|\mathcal S'|=p+1>\frac{d}{2}   \ge   \sqrt{(n-1)} \ge \sqrt{m}$, by (\ref{eqvertices}). 
\qed

\vspace{3mm}
  We now show how to extract free subgroups of rank  $r\geq m'$ from any disjoint pS-family $\mathcal S' =\{ h_1,...,h_{m'}\}$ whose axes are in one of the configurations (c) or (d).

\begin{lem}\label{pS-disjoint_smallproj}
	Let $\mathcal S'=\{h_1,..., h_{m'}\}$ be a disjoint pS-family with small projections and such that $\tau(h_i)=\tau$ for all $i$. Then $\mathcal S'$ generates a Schottky  group.
\end{lem}

{\em Proof}.
     We set $X_i=\{\mathsf v\in\mathcal T\,|\, p_i(\mathsf v)\not\in I_i\}$ for any $i=1,.., m'$ . By assumption $p_i(\ax(h_k))\in I_i$ for any $k\neq i$ and thus we conclude that $X_i\cap X_k=\varnothing$ (otherwise $\mathcal T$ would contain a loop).  Now observe that if $\mathsf v\in X\smallsetminus X_i$, then $p_i(\mathsf v)$ is in  $I_i$. Since $\tau(h_i)=\tau>\diam(I_i)$ and $p_i(h_i(\mathsf v))=h_i(p_i(\mathsf v))$  we conclude that $h_i^{\pm 1}(\mathsf v)\in X_i$. Thus $\langle h_1,.., h_{m'}\rangle$ is a Schottky subgroup.
\qed

\begin{lem}\label{pS-disjoin_path}
	Let $\mathcal S'=\{h_1,...,h_{m'}\}$ be a disjoint sequential  pS-family such that $\tau(h_i)=\tau$ for all $i$. 
	Then there exist  $\varepsilon_2, ..., \varepsilon_{m'-1} \in \{ \pm1\}$, such that the family 
 $$\mathcal S''=\{ h_1,...,h_{k}^{\varepsilon_{k}}h_{k+1}h_{k}^{-\varepsilon_{k}},..., h_{m'}\}$$
for $k=2,...,m'-1$,   generates a Schottky subgroup.
\end{lem}

{\em Proof}.
    For each $i=2,.., m'-1$ choose $\varepsilon_i$ so that the segment $[\mathsf n_{i}(i-1), \mathsf n_{i}(i)]$ is coherently oriented with the translation direction of $h_i^{\varepsilon_i}$.
    Now  consider the set  $\mathcal S''=\{ h_1,...,h_{k}^{\varepsilon_{k}}h_{k+1}h_{k}^{-\varepsilon_{k}},..., h_{m'}\}$, 
    and notice  that it has a single nerve which is  given by 
    \small
    $$[ \, \mathsf n_{1}(1)\, ,\,  h_2^{\varepsilon_2}\mathsf n_3(2) \, ]  \;  
        \bigcup_{k=2}^{m'-1} [  \, \mathsf n_{k}(k-1) \,,\,  h_k^{\varepsilon_k} \mathsf n_{k+1}(k) \,]  \;\, \bigcup \;\,
        [ \, h_{m'-1}^{\varepsilon_{m'-1}}   \mathsf n_{m'}(m'-1) \,,\,  \mathsf n_{m'}(m'-1) \,]   $$
    \normalsize 
 \noindent (see Fig. \ref{schottky_from_disjoint_paths}). 
    We can then apply Lemma \ref{pS-disjoint_smallproj} to conclude that the group generated by $\mathcal S''$ is a Schottky  group.
\qed

\bigskip
\subsection{Proof of the Entropy-cardinality inequality, Theorem \ref{ent-card_ineq}.}
\label{secproofEC} ${}$ 
\vspace{1mm}
 
\noindent By Proposition \ref{unif_bounded_large_pS-family} , for any generating set $S$ with $|S|=n$ we can  find a pS-family $\mathcal S$ of cardinality $m \geq \sqrt[16]{n}/\sqrt 2$ 
such that $\tau(h)=\tau$  and   $|h|_S\le 20\,(k+3)$  for all $h \in \mathcal S$. 
Moreover, by the same theorem, either there exists $h_1\in\mathcal S$ such that $\ax(h_i)\cap\ax(h_1)\neq\varnothing$ for all $i$,  or the axes of elements in $\mathcal S$ are pairwise disjoint.  \\
If  the axes of all the elements in $\mathcal S$ intersect $\ax(h_1)$, we know by Proposition \ref{extraction_intersecting} that there exists a pS-subfamily $\mathcal S'\subset\mathcal S$  which  either is  agglomerated with cardinality $m' \ge\sqrt{m}$,  or is sparse with respect to $h_1$ and has cardinality $m' \ge\sqrt{m}+1$.
In both cases, we deduce from Lemmas \ref{pS_bigcap} and \ref{pS-almostdisjoint} a free Schottky subgroup $H$ 
  of rank $  m'  \geq \sqrt[32]{n}/\sqrt[4]2 $,
 generated by   elements having $\mathcal S'$-length at most $4$.  \\
On the other hand, if   the axes of  $\mathcal S$ are pairwise disjoint, by Proposition \ref{disjoint_schottky} there exists a disjoint, pS-subfamily $\mathcal S' \subset \mathcal S$ of cardinality $m' \ge \sqrt{m}$ which is either   with small projections  or   sequential.
We then deduce again,  from Lemmas  \ref{pS-disjoint_smallproj} and \ref{pS-disjoin_path},   the existence of a  Schottky subgroup $H$ 
of rank 
$ m' \geq \sqrt{m} \geq \sqrt[32]{n} / \sqrt[4]{2} $,
generated by elements having $\mathcal S'$-length at most $3$.
Hence,  
$$\ent(G,S) \ge \ent(H, d_S) \ge  \frac{1}{80(k+3)}  \ent(H, \mathcal S')
                                                   \ge  \frac{1}{80(k+3)} \log \left(2m'-1\right)$$                                                 
 which gives the announced inequality.  \qed\\

\section{Applications} \label{ssapplications}


%
 

\vspace{1mm}
Recall that, given a marked group $(G,S)$, a complete set of relators for $G$ is  a finite subset $R$ of the free group on $S$ such that $G \cong \mathbb F(S)/\langle\!\langle R\rangle\!\rangle$, where $\langle\!\langle\,R\,\rangle\!\rangle$ denotes the normal subgroup generated by the elements in $R$.

A first  consequence of the Entropy-Cardinality inequality is a general, finiteness theorem for abstract groups admitting acylindrical splittings, which has an interest in its own:

\begin{thm}\label{finiteness_short_relators} 
	The number of  isomorphism classes of marked groups  $(G,S)$ \linebreak  admitting a non-elementary, $k$-acylindrical splitting and a complete set of relators of length less than $\ell$,  satisfying $\ent(G,S) \leq E$ is finite, bounded by  an explicit function $N(k,\ell, E)$.
\end{thm}

{\em Proof of Theorem \ref{finiteness_short_relators}}.
The Entropy-Cardinality inequality yields   a corresponding  bound on the cardinality of $S$ in terms of the acylindricity constant and of $E$:
	\small
	$$|S| \le  \left( 1+ e^{E\cdot(80k+3)} \right)^{32}  =n(k,E)$$
	\normalsize
Now, the number of possible presentations by relators of length smaller than $\ell$ on an  alphabet ${\mathcal A}$ of at most $n$ letters is roughly bounded by $2^{n+n^\ell}$ (that is the number of subsets $S$ of ${\mathcal A}$, times the number of subsets $R$ of ${\mathcal A}^\ell$); this  gives the announced bound of the number of marked groups $(G,S)$ by the function $N(k,\ell, E) = 2^{n(k,E)+n(k,E)^\ell} $.
\qed

   \vspace{3mm}
  In the following, we will make use of these basic facts about amalgamated products and HNN-extensions over malnormal subgroups (the proof of Lemma \ref{malnormalAP}(i) and (ii)  can be found,  with the original terminology,  in \cite{karso}, Section \S5; the case of HNN-extensions is analogous, cp. \cite{cer_pi1} for more details):

\begin{lem}\label{malnormalAP}
	Let $A \, \ast_C B$ be the amalgamated product of two groups $A,B$. \\
	Assume that  $C$  is a proper subgroup having  index greater than $2$ in $A$ or in $B$,  and let  $\iota_A: C \f A$ and $\iota_B: C \f B$ be the natural inclusions: 
	 
\begin{itemize}[leftmargin=8mm]
\item[(i)]  if one between $\iota_A(C)$,  $\iota_B(C)$ is malnormal in the respective group, then   $A   \ast_C B$ \linebreak is a $1$-step malnormal amalgamated product,  and has a non-elementary \linebreak $2$-acylindrical splitting;  
\item[(ii)]  if both $\iota_A(C)$ and $\iota_B(C)$ are malnormal in the respective groups, then  $A   \ast_C B$ \linebreak is a $0$-step malnormal amalgamated product, and has  a non-elementary \linebreak $1$-acylindrical splitting.
\end{itemize}
 \end{lem}


\begin{lem}\label{malnormalHNNext}
	Let $A \, \ast_\varphi =\langle A, t \;   |\; t^{-1}ct=\varphi(c) \rangle$ 
be the HNN-extension of $A$ with respect to an isomorphism $\varphi: C_{-} \f C_+$ between two subgroups  one of which has index at least two in $A$, and let  $\iota: C_- \f  \, A \ast_\varphi $ be the natural inclusion: 
\begin{itemize}[leftmargin=8mm]
\item[(i)] if $\iota (C_-)$ is malnormal in $A\, \ast_\varphi$, then $C_+$ and $ C_-$ are malnormal and conjugately separated in $A$;
\item[(ii)] if  the subgroups $C_+, C_-$ are malnormal and conjugately separated in $A$, then   $A\, \ast_\varphi$ is a 
 $0$-malnormal HNN-extension, hence it, admits a non-elementary \linebreak $1$-acylindrical splitting.\\
%

\end{itemize}
\end{lem}
 

\vspace{-2mm}
\subsection{Quasiconvex groups of $\delta$-hyperbolic and  $CAT(0)$-spaces}\label{ssdeltahyp}${}$

\noindent Recall that an action of a group $G$ on a $\delta$-hyperbolic or $CAT(0)$-space $X$ is called {\em quasi-convex} if there exists an orbit $S=Gx_0$ which is a $D$-quasi-convex subset of $X$ (i.e. all the geodesics joining  two points $x_1,x_2 \in S$ are  included in the closed $D$-neighborhood of $S$), for some $D>0$. One also says that $G$ is   a {\em quasi-convex group of $X$}.

\vspace{3mm}
  {\em Proof of Corollary \ref{deltahyp}}.
\noindent  Any marked $\delta$-hyperbolic group $(G,S)$ possesses a complete set of relators of length $\ell\le (4\delta+6)$ (see for instance \cite{BH}, Chapter III.$\Gamma$, Proposition 2.2):  the conclusion  then immediately follows from Theorem \ref{finiteness_short_relators}.\qed
\vspace{3mm}

   {\em Proof of Corollary \ref{quasiconvex}}.
Let $G$ act  on  a  proper, geodesic space $(X,d)$,  and let $x_0 \! \in \!X$ a point with $D$-quasi-convex orbit.
By a classical argument, the set \linebreak $S =  \{g\in G\,|\, d(x_0, g.x_0)\le 2D+1\}$ generates $G$. 
 Actually, for every $g\in G$, consider  a geodesic $c: [0,\ell] \rightarrow X$ from $x_0$ to $g.x_0$, and orbit points $g(k). x_0$ such that $g(0)=1$, $g(\lceil \ell \rceil)=g$ and  $d(c(k), g(k). x_0) \leq D$, 
  given by the condition of $D$-quasiconvexity.  Then, setting  $\gamma(1)= g(1)$ and  
  $\gamma(k) = g(k-1)^{-1}g(k)$, one has that the $\gamma(k)$'s are in $S$ and 
 $g=\gamma(1) \cdot \gamma(2) \cdots \gamma(\lceil \ell \rceil)$.
By construction, we also have
\begin{equation}\label{eqquasiisom}
\frac{1}{2D +1}  d(x_0, g\cdot x_0) \leq |g|_S \leq  d(x_0, g. x_0) +1
\end{equation}
therefore the marked group $(G,S)$ is  $(2D+1, 1)$-quasi-isometric to  the orbit $G.x_0$, endowed with the distance $d$ induced by $X$;  it follows that, in case (i), $(G,S)$ is   $\delta'$-hyperbolic,   for some $\delta'=\delta'(\delta, D)$. On the other hand,  from the left-hand side of (\ref{eqquasiisom}) we deduce that  
$$\ent(G,S)\le (2D+1) \cdot  \ent(G \curvearrowright X)\le (2D+1)E$$ 
The conclusion in case (i) then follows from Corollary \ref{deltahyp}.\\
In case (ii), we proceed similarly to \cite{BH} Chapter III.$\Gamma$, Proposition 2.2, by replacing the cocompactness assumption by quasiconvexity. Namely, given $S$ as above, one decomposes any relation $r=s_1\cdots s_n$  on $S$ for $G$ as a product of relations  
$$r = \prod_{i=1}^{n-1} \left( \sigma_i s_{i+1} \sigma_{i+1}^{-1} \right) $$
where  $\sigma_i:=s_1\cdots s_i$, and then shows  that each relation $\sigma_i s_{i+1} \sigma_{i+1}^{-1}$ is product of conjugates of words $r_k$ of $S$-lenght at most $8D+6$. Actually,  choose again geodesics $c_i:  [0,\ell_i] \rightarrow X$  from $x_0$ to $\sigma_i . x_0$, and then orbit points $g_i(k) . x_0$ with   $g_i(0)  =1$, $g_i(\lceil \ell_i \rceil) = \sigma_i$  for all $i$ and  $d(c_i(k) , g_i(k) . x_0) \leq D$, provided by the $D$-quasiconvexity. We then consider the elements  $\gamma_i(k) := g_i(k)^{-1} g_i(k+1) $ and $\mu_i(k) := g_i(k)^{-1} g_{i+1}(k) $ of $G$, and notice that  
$$d(x_0, \gamma_i(k) . x_0) 
\leq 2D + d(c_i(k) , c_{i}(k+1) ) = 2D +1 \,.$$
and that, by the convexity of the metric of $X$,  
$$d(x_0,\mu_i(k) . x_0)  
\leq 2D + d(c_i(k) , c_{i+1}(k)  )  \leq    2D + d(\sigma_i . x_0, \sigma_{i+1}  . x_0)  =  4D +1$$
Therefore,   $\gamma_i(k) \in S$ and $|\mu_i(k)|_S \leq 4D+2$,   by (\ref{eqquasiisom}); so, $\mu_i(k)$ can be represented by a word $\tilde \mu_i(k)$ on $S$ of length  $\leq 4D+2$.
Accordingly,  all  the relations $\sigma_i s_{i+1} \sigma_{i+1}^{-1}$ can   be decomposed as  products of conjugates
$$\sigma_i s_{i+1} \sigma_{i+1}^{-1}=   g_i(\lceil \ell_i \rceil) r_i (\lceil \ell_i \rceil)  g_i(\lceil \ell_i \rceil) ^{-1}  \cdots g_i(k) r_i (k)   g_i(k) ^{-1} \cdots g_i(1) r_{i} (1)  g_i(1) ^{-1}$$
where the words $r_i (k) :=  \tilde \mu_i(k+1)   \gamma_{i+1}  (k)^{-1}   \tilde \mu_i(k)^{-1}    \gamma_{i}  (k) $ represent relators    on $S$  whose  $S$-lengths do  not exceed  $8D+6$.
The conclusion for case (ii) then follows from Theorem \ref{finiteness_short_relators}.\qed

%

\vspace{5mm}
\subsection{CAT(0)-spaces with negatively curved splittings}\label{ssnegcurvsplittings}${}$ 

 
\noindent We say that a locally $CAT(0)$-space $X$  {\em admits a splitting} if $X$ is isometric to the gluing $Y_1 \sqcup_{\phi} Y_2$ of two locally $CAT(0)$-spaces $Y_1, Y_2$ along  compact, locally convex, isometric  subspaces $Z_i \hookrightarrow X_i$ via an isometry $\phi: Z_1 \rightarrow Z_2$; or if $X$ is isometric  to the space $Y \sqcup_{\phi} $ obtained by identifying two such  subspaces $Z_i \subset Y$ to each other by an isometry $\phi$. 
The splitting is {\em non-trivial} if the corresponding splitting of $\pi_1(X)$ as an amalgamated product or HNN-extension is non-trivial. \linebreak
Notice that the  space obtained by such gluings is always  locally $CAT(0)$ (cp. \cite{BH}, Prop. 11.6). 
We will say that $X$ has a {\em negatively curved splitting} if  the subspace $Z$ of $X$ obtained by identifying $Z_1$ to $Z_2$ has a neighbourhood  $U(Z)$ in $X$  such that $U(Z)\setminus Z$   is a locally $CAT(-k)$-space for some $k>0$.
 \vspace{1mm}

The following fact  is crucial to prove acylindricity for negatively curved splittings of $CAT(0)$-spaces, and we believe it is folklore; we will give a proof of it in Appendix \ref{appcat(0)}, by completeness.
  \vspace{2mm}
 
\noindent {\bf Proposition  \ref{malnormal}}
{\em
	Let $Z$ be a  compact,   locally convex subspace of a compact, complete locally  $CAT(0)$-space $X$.
	Assume that $X$ is negatively curved around $Z$: then, $H=\pi_1 (Z)$ is malnormal in   $G=\pi_1 (X)$.
}



  \vspace{4mm}
  {\em Proof of Theorem \ref{thmcat(0)}}.
If $X$ admits 	a non-trivial, negatively curved splitting \linebreak  $X=X_1 \sqcup_{Z} X_2$  or $X=X_0 \sqcup_{\phi} $, along two  isometric, locally convex  subspaces $Z_i \cong Z$  identified   to each other via an isometry $\phi$, 
 then $\pi_1 (X)$ splits as a non-trivial amalgamated product of the  groups $G_i = \pi_1 (X_i)$ over $H=\pi_1(Z)$, or as a non-trivial HNN-extension of $G_0=\pi_1 (X_0)$   
along  subgroups $H_i=\pi_1 (Z_i)$, via the isomorphism $\phi_\ast: H_1 \f H_2$.
 The subgroups $H_i$ and $H$,  are malnormal respectively in each  $G_i$ and in $G$,  
 by Proposition \ref{malnormal}.
Moreover,  each $H_i$ does not have index $2$  in $G_i$, 
or it would be normal and malnormal, hence trivial (since the splitting is supposed to be non-trivial).
Therefore,   $\pi_1 (X)$ admits a non-elementary $1$-acylindrical splitting.
By Corollary \ref{quasiconvex} , we deduce that  $\pi_1 (X)$ belongs to a finite class of groups, whose number is bounded by a universal function of $E$ and $D$. 
Since the locally $CAT(0)$-spaces are aspherical, we can conclude by Whitehead's theorem the finiteness of the homotopy types.\qed
\vspace{3mm}

  {\em Proof of Corollary \ref{cormanifoldsk<0}}.
It follows from the fact that,   in dimension greater than $4$, the homeomorphism type of closed, non-positively curved manifolds is determined by their homotopy type, from the solution of the Borel Conjecture by Bartels-L\"uck  for $CAT(0)$-manifolds  \cite{bartels-luck}. Moreover,  the works of Kirby-Siebenmann  \cite{kirby-siebenmann} and of Hirsch-Mazur \cite{hirsch-mazur}
 on PL structures and their smoothings imply the finiteness  of smooth structures in  dimension $n\geq 5$ (cp. also 
 \cite{lance}, Thm. 7.2). \linebreak
 The fact that the fundamental group determines the  diffeomorphism type is well-known in   dimension  $2$,  and follows in dimension $3$ from Perelman's solution of the geometrization conjecture (we now know that any closed, negatively curved $3$-manifold  also admits  a hyperbolic metric; so Mostow's rigidity applies).
 \qed

 
    

  \vspace{2mm}
\subsection{Two-dimensional orbifolds}${}$ 
\label{ss2dimorbifolds}

\noindent We recall shortly some basic facts about orbifolds;  for a primer concerning $2$-dimensional orbifolds we refer to \cite{Scott} and \cite{Thurst}.  
Following Thurston \cite{Thurst},  a $n$-dimensional orbifold $\mathcal O$ (without boundary) is a Hausdorff, paracompact space which is locally homeomorphic either  to   $\mathbb{R}^n$, or to the  quotient of $\mathbb{R}^n$  by a finite group action;  similarly,  $n$-orbifolds with boundary   also have points whose neighbourhood is  homeomorphic to  the quotient of the half-space $\R_+^n$ by a finite group action.
 For the sake of simplicity we shall consider uniquely compact $2$-dimensional orbifolds with {\em conical singular points}, that is points which have a neighborhood modelled on the quotient of $D^2$ by a finite {\em cyclic} group. 
 Nevertheless, it follows from the description of singularities in \cite{Scott} and \cite{Thurst} that given a general compact $2$-orbifold there exists a canonically constructed double cover which has only conical singularities; this cover is obtained by doubling the underlying space along the {\em reflector lines}, duplicating the conical singular points and trasforming the so called {\em corner reflectors} into conical singular points.

We shall denote by $\mathcal O=\mathcal O(g,h;p_1,...,p_k)$ the compact $2$-orbifold having as underlying topological space a compact surface $|\mathcal O|$ of genus $g\in \Z$ (using  negative values for the genus of non-orientable surfaces),  $h$ boundary components  and $k$ singular points $x_1,.., x_k$ of orders $p_1,...,p_k$. 
By the classification of  compact $2$-orbifolds given by Thurston in \cite{Thurst}, Ch.~13, an orbifold with conical singular points is completely determined by its underlying topological space together with the number and the orders of its singular points. In view of this fact, we shall say that two smooth compact $2$-orbifolds with conical singularities are {\em isomorphic} if they have the same underlying surface, the same number of singular points, and the same order at each singular point, up to permutations. 

 For a formal definition of the orbifold fundamental group, we refer to \cite{BMP}; to our purposes, it will be sufficient to recall that the orbifold fundamental group of   $\mathcal O=\mathcal O(g,h;p_1,...,p_k)$   admits one of the following presentations: \pagebreak
\small
$$\pi_1^{orb}(\mathcal O)=\left\langle  a_1, b_1,..,  b_g,  c_1,.., c_k,  d_1,.., d_h\,\left |\,\begin{array}{l}
\prod_{i=1}^{g}[ a_i, b_i]\cdot\prod_{j=1}^{k}  c_j\cdot\prod_{\ell=1}^{h}  d_\ell=1\\ 
 c_i^{p_i}=1
\end{array}
\right.\right\rangle$$
$$\pi_1^{orb}(\mathcal O)=\left\langle  a_1,...\, a_{|g|},   c_1,..., c_k,  d_1,..., d_h\,\left |\,\begin{array}{l}
\prod_{i=1}^{|g|} a_i^2\cdot\prod_{j=1}^{k}  c_j\cdot\prod_{\ell=1}^{h}  d_\ell=1\\ 
 c_i^{p_i}=1
\end{array}
\right.\right\rangle$$
\normalsize
depending on whether the genus is positive or negative. The generators $a_i$'s, $b_j$'s, $d_\ell$'s are the fundamental system of generators of  $|\mathcal O|$, with $d_\ell$ corresponding to the $l$-th boundary compoent,  whereas the $c_i$'s represent  the generators of the isotropy groups associated to  the singular points of $\mathcal O$.

\noindent The usual Euler characteristic can be generalized in a natural way to the case of compact $2$-orbifolds (\cite{Scott}), and in the case  we are considering the formula reads:  
\small
$$\chi_{orb}(\mathcal O(g, h; p_1,..., p_k))=\chi(|\mathcal O|)-\sum_{i=1}^k\left(1-\frac{1}{p_i}\right)=\left\{\begin{array}{c}
2-2g-h-\sum_{i=1}^k\left(1-\frac{1}{p_i}\right) \\
2-|g|-h-\sum_{i=1}^k\left(1-\frac{1}{p_i}\right) 
\end{array}\right.$$
\normalsize
\noindent depending on the sign of the genus $g$.
The orbifold Euler characteristic is useful to distinguish between those orbifolds which are finitely covered by a compact surface (usually referred as {\it good orbifolds}) and those who are not (the {\it bad} ones): it follows from \cite{Thurst} (Theorem 13.~3.~6)  that if $\mathcal O(g, h; p_1,..., p_k)$ has positive orbifold Euler characteristic and  $LCM(p_1,..., p_k)\cdot\chi_{orb}(\mathcal O) \neq 1,2$ then it is a {\it bad orbifold}. Moreover, the Euler characteristic a good orbifold $\mathcal O$  determines the kind of {\em geometric structure} that can be given to $\mathcal O$:  that is,  whether  its universal cover $\widetilde{\mathcal O}$ --which is a simply connected surface-- admits a $\pi_1^{orb}(\mathcal O)$-invariant spherical, flat or hyperbolic Riemannian metric. 
We shall accordingly call compact $2$-orbifolds with negative Euler characteristic {\it orbifolds of hyperbolic type}. 

Among compact $2$-orbifolds of  hyperbolic type there is a particular  family whose  groups are generated by symmetries of hyperbolic triangles along their edges: these groups are called  {\it  triangle groups}, and the corresponding orbifolds are  those of the form $\mathcal O(0,0;p,q,r)$, with $p=2$, $q =3$, $r \geq 5$ or with 
$p,q,r \geq 3$ and at least one of them strictly greater than $3$.
\vspace{2mm}

%

In Theorem  \ref{fntn_orb} we  consider general {\em Riemannian}   $2$-orbifolds: that is,  good  compact $2$-orbifolds $\mathcal O$ (with conical singularities) whose orbifold universal cover $\widetilde{\mathcal O}$  is endowed with {\em any} $\pi_1^{orb}(\mathcal O)$-invariant Riemannian metric (that is, not necessarily with constant curvature).
The entropy of $\mathcal O$ is correspondingly defined, as explained in the introduction, as $\ent (\mathcal O) = \ent (\pi_1^{orb}(\mathcal O) \curvearrowright \widetilde{\mathcal O})$.

Notice that  compact, two-dimensional orbifolds with zero orbifold  Euler  characteristic   yield only finitely many isomorphism classes (cp. \cite{Thurst}, Theorem 13.3.6),  as it can be checked directly from the Euler characteristic formula,   without any assumption on their entropy and diameter.

The proof Theorem  \ref{fntn_orb} will then split in two separate cases:\\
(a)  non-triangular, Riemannian $2$-orbifolds of hyperbolic type; \\
(b)  triangular $2$-orbifolds of hyperbolic type.\\
 In fact, in the first case, the orbifold groups admit a $2$-acylindrical splitting, as we will show in detail in Proposition \ref{splitting2orbifolds} of Appendix \ref{apporbifold}, and we can use the Entropy-cardinality inequality. In the second case,   acylindrical splittings are not available for such groups,  and we will need an ad-hoc  computation  to conclude.
 We will denote  by $\mathscr O_{h,nt}^2(E,D)$ and  $\mathscr O_{h,t}^2(E,D)$ the classes  of compact,  Riemannian $2$-orbifolds of negative orbifold Euler characteristic with  entropy and diameter  bounded by $E$ and $D$,   which fall, respectively, in cases (a) and (b).
\vspace{4mm}

   {\em Proof of finiteness of isomorphism types in $\mathscr O_{h,nt}^2(E,D)$.} \\
 The orbifold groups in this class  admit a $2$-acylindrical splitting, cf. Proposition \ref{splitting2orbifolds} of Appendix \ref{apporbifold}.
Then, for any orbifold $\mathscr O \in  \mathscr O_{h,nt}^2(E,D)$ we choose $\tilde x_0 \in \widetilde{\mathscr O}$, 
and apply  to $G=\pi_1^{orb} (\mathscr O)$ and to the open ball $U=B_{\widetilde{\mathscr O}}(\tilde x_0, 2D+\epsilon)$ of  $\widetilde{\mathscr O}$, for arbitrary $\epsilon>0$,   the  following, classical result (see, for instance,  \cite{Serre}, p.30):

	\begin{lem}\label{SPorbifold}
		Let $G$ act by homeomorphisms on a path-connected, simply connected topological  space $X$, and  let $U$ be a path-connected open set such that $G.U=X$.\linebreak
		Let $S=\{s\in G\,|\, sU\cap U\neq\varnothing\}$ and $T = \{ (s_1,s_2) \in S \times S\;|\; U \cap (s_1U) \cap (s_1s_2 U) \neq \emptyset \}$. \linebreak
 Then  $S$ generates $G$, and $G \cong F(\Sigma) \slash \langle \! \langle \Theta \rangle \! \rangle$
		where  $\Sigma$ is the set of symbols $ \{ x_s  \,|\, s \in S\}$,  \linebreak and $\Theta$ is the set of words on $S\cup S^{-1}$ given by $\{ x_{s_1} x_{s_2} x_{s_1s_2}^{-1} \,|\, (s_1,s_2) \in T\}$.\\
 (Notice that if $(s_1,s_2) \in T$, then  $s_1s_2 \in S$ so $x_{s_1s_2}^{-1}$ makes sense). 
		\end{lem}

\noindent  For $M =2D+\epsilon$,  the lemma yields a finite  generating set $S_M$ of $G$, such that  $d(x_0, g \cdot x_0) < M$ for all $g \in S_M$, 
 which we call {\em $M$-short generators of $G$ at $x_0$}, and a {\em triangular presentation} of the  group: that is, such that the group of relations is generated, as a normal subgroup of  $ \mathbb{F} (S_M)$, by relators of length at most $3$.\\
 Since $d(x_0, g\cdot x_0) \leq M \cdot |g|_{S_M}$, we have  $\ent(\pi_1(\mathscr O), S_M) \leq M \ent(\mathscr O) \leq ME $.
Letting $\epsilon \rightarrow 0$, we deduce from Theorem \ref{ent-card_ineq}  that
  $|S_M|  \le \left(e^{800\,DE} +1\right)^{32} =N(E,D)$.
  \\
As in the proof of Theorem \ref{finiteness_short_relators} observe now that the number of possible triangular presentations that can be build with letters from some subset $S$ of  an alphabet of cardinality  $N$   does not exceed $2^{N+N^3}$;  therefore, the number of  such orbifold groups  is bounded from above by $2^{N(E,D)}+2^{N(E,D)^3}$. 
 To conclude,  remark that    $\pi_1^{orb} (\mathscr O)$  determines the isomorphism class of $\mathscr O$ for closed orbifolds, while for orbifolds with boundary the isomorphism class is determined by the orbifold group and the number of boundary components; in any case, there are a finite number of non-isomorphic $2$-orbifolds  for any given group $\pi_1^{orb} (\mathscr O)$.\qed
  \vspace{4mm}
    

  {\em Proof of finiteness of isomorphism types in $\mathscr O_{h,t}^2(E,D)$.} \\     
In this case, we call $G_{p,q,r}=\pi_1^{orb}( \mathcal O(0,0;p,q,r))$ the triangle orbifold group and consider the fixed point $\tilde x_0\in \widetilde{\mathcal O}(0,0;p,q,r)$  of the torsion element $c$ of order $ r\geq4$.
To evaluate   the Poincar\'e series of $G_{p,q,r}$  at $\tilde x_0$, we need the following

\begin{lem}\label{injectivity_Wr}
	Let $G_{p, q, r}=\langle a, b, c\,|\, abc=a^{p}=b^q=c^r=1 \rangle$ where $3\le p\le q\le r$ and $r\ge 4$ be the fundamental group of a triangular orbifold of hyperbolic type. \linebreak The following set of elements of $\Z_{p}*\Z_{r}$ naturally injects   into $G_{p, q, r}$:
	\small
	\begin{equation*}
	\mathcal W_{r}=\left\{ c^{j_1}a c^{j_2}a\cdots c^{j_m}a\,|\, m\in\N,\, 2\le j_i\le\frac{r}{2}\right\}
\end{equation*}
	\normalsize
(By a slight abuse of notation, we will use  $\mathcal W_r$ and the reduced forms $c^{j_1}a\cdots c^{j_m}a$ to denote both the elements of $Z_p \ast Z_r$ and their images in $G_{p, q, r}$).
	
\end{lem}

\noindent Assuming this lemma for a moment, we can compute:
\begin{equation}\label{eqpoincare}
\quad P_s (G_{p,q,r}, \tilde x_0)
\ge\sum_{g\in\mathcal W_r} e^{-s\, d(\tilde x_0, g. \tilde x_0)}
=\sum_{n }\sum_{2\le j_i\le \frac{r}{2}} e^{-s \, d(\tilde x_0,   c^{j_1} a\cdots c^{j_n} a . \tilde x_0)}
\end{equation}
\noindent Moreover, since we are considering triangular orbifolds of diameter bounded by $D$ and $\tilde x_0$ is the fixed point of $c$, we have:
\vspace{-3mm}

$$ d(\tilde x_0,  c^{j_1}a\cdots c^{j_n}a . \tilde x_0)
\le  n d(o, a. o) + \sum_{i=1}^n d(o, c^{j_i}. o)
\le 2nD$$

\vspace{-1mm}
\noindent which plugged   into (\ref{eqpoincare}) yields:
\small
    $$P_s (G_{p,q,r}, \tilde x_0) \ge \sum_{n }\sum_{2\le j_i\le\frac{r}{2}} e^{-2sn\,D}\ge \sum_{n>0}\left(\left[\frac{r}{2} -1\right] e^{-s\,2D}\right)^n$$
    \normalsize
Since $P_s (G_{p,q,r}, \tilde x_0)$ converges for all $s >E$, this shows that  $[\frac{r}{2} -1] \leq e^{ 2ED}$. \linebreak As $r\ge q\ge p$ this proves the finiteness of the  fundamental groups and of isomorphisms classes of triangular orbifolds with bounded entropy and diameter.\qed
 \vspace{4mm}
 
  {\em Proof of Lemma \ref{injectivity_Wr}}.
By the canonical presentations of the compact, $2$-orbifold groups recalled before, we know that  $G_{p, q, r}=\Z_{p}*\Z_r/ \langle\!\langle (ca)^q\rangle\!\rangle$, since $b^{-1}=ca$.  \linebreak
We start with the case  where $q \geq 4$. 
Following \cite{lynschu},  we call a word  on \linebreak $S=\{ a,a^{-1}, c, c^{-1} \}$ with normal form $w = a^{p_1}c^{q_1}\cdots a^{p_n}c^{q_n}$, possibly with $p_1 =0$ or $q_n=0$, {\em weakly cyclically reduced} if the last syllable of $w$ (that is, $c^{q_n}$, if $q_n \neq 0$) is different from the first one.
If $R$ denotes the set of weakly cyclically reduced conjugates of $(ca)^{\pm q}$, then the quotient   $G_{p, q, r} = \Z_{p}*\Z_r/ \langle\!\langle R \rangle\!\rangle$  satisfies  the $C'(\frac{1}{6})$ condition for free products 
(every prefix of an element $ r \in R$ which is a piece has syllable length
\footnote{Notice that the small cancellation theory on free products differs from general cancellation theory, the relevant length and notion of {\em piece} being given by the syllable length and  by the subdivision in syllables provided by the normal form.}
 smaller than $\frac{1}{6} \ell (r)$).
 It follows from small cancellation theory  that any element $w\in \Z_p*\Z_r$  belonging to $\langle\!\langle R\rangle\!\rangle$ has a normal form which contains as a subword a prefix $r_0$ of some element in $R$ of syllable length $\ell(r_0) \geq 5$. 
 It is  straightforward to check that none of the elements of $\mathcal W_r \mathcal W_r^{-1}$ contains such a subword, hence  $\mathcal W_r$ injects  into $G_{p,q,r}$. \\
 The argument when  $p=q=3$ and $r \geq 4$ is the similar, with the difference  that in this case the set $R$ of weakly cyclically reduced conjugates of $(ca)^{\pm q}$ does not satisfy the $C'(\frac{1}{6})$ condition, but   conditions $C'(\frac{1}{4})$ and  $T(4)$ (given $r_1, r_2, r_3\in R$, at least one among  $r_1r_2$, $r_2r_3$, $r_3r_1$  is {\em semi}-reduced, cf. \cite{lynschu}).
 Nevertheless, if $R$ satisfies $C'(\frac{1}{4})$ and $T(4)$, it is still true that any reduced word $w\in\langle\!\langle R\rangle\!\rangle$ contains a prefix $r_0$ of an element in $R$ as subword,   with syllable length $\ell(r_0) \geq 4$ (\cite{lynschu}, Ch.5,  Thm.4.4). We then proceed as before and deduce that $\mathcal W_r$  injects into $G_{p,q,r}$.\qed


  \vspace{4mm}
\subsection{Non-geometric $3$-manifolds} ${}$ 
\label{sec3manifolds}

\noindent  In this section we will prove the finiteness results Theorem \ref{glob_fnt} and Corollary \ref{fnt_HT} for non-geometric Riemannian $3$-manifolds  with bounded entropy and diameter. \linebreak
We will first show the finiteness of fundamental groups of manifolds in the class $\mathscr M_{ng}^{\partial}(E,D)$,  and then  explain how to deduce Corollary \ref{fnt_HT} from Theorem \ref{glob_fnt}.  

\vspace{1mm}
The proof of the finiteness of the isomorphism classes of  fundamental groups in  $\mathscr M_{ng}^{\partial}(E,D)$ relies on the fact that  the fundamental group of every  non-geometric $3$-manifold,  closed or compact with non-spherical boundary components, admits  a non-elementary $4$-acylindrical splitting: this fact was proved by \cite{WiZa} (see also \cite{cer_pi1} for further details).
 The splitting is relative  either to  the decomposition of $\pi_1(X)$ as a free product given by prime decomposition, or to  the decomposition  of  $\pi_1(X)$  as an amalgamated product over rank 2, abelian subgroups  provided  by the   JSJ-decomposition, for irreducible manifolds, cf. \cite{CS}, Section \S4.
We can then  apply the  Entropy-Cardinality inequality to the classical triangular presentation of $\pi_1(X)$ given by the Lemma \ref{SPorbifold}, and proceed as in the proof of Theorem  \ref{fntn_orb}, case (a).


\vspace{3mm}
{\em Proof of Theorem \ref{glob_fnt}.}
	For any $X \in \mathscr M_{ng}^{\partial}(E,D)$  pick  $x \in X$ and let $M=2D+\epsilon$. 
	Then, consider   the set $S_{M}$   of   $M$-short generators  at  $x$.
	 As $\pi_1(X)$ has a $4$-acylindrical splitting, it follows from   from Theorem \ref{ent-card_ineq}  that
  \vspace{-2mm}
  \small
  
  	$$E \ge  \ent(X)\ge\frac{1}{M}\ent(\pi_1(X), S_M )
	\ge \frac{1}{560M}\log\left(\sqrt[32]{|S_M|}-1\right)$$
  
  \normalsize
  \vspace{-1mm}
  
  \noindent   	and letting $\epsilon \rightarrow 0$ we obtain  $|S_M| \le \left(e^{1120\,E\,D} +1\right)^{32} =N(E,D)$.
	Therefore, $X$ admits a triangular presentation on a generating set of cardinality at most $N(E,D)$.
	Since the number of possible triangular presentations that can be build with letters from some subset $S$ of an alphabet ${\mathcal A}$ of $N$ letters does not exceed $2^{N+N^3}$, this concludes the proof.
\qed

\vspace{2mm}
Now, the following statement is  consequence of several results of $3$-dimensional geometry and topology. 
Since it relies on facts which are now folklore (and are sometimes only sketched in literature), we will provide  a full proof  in Appendix \ref{appJS}, together with all the references and the $3$-dimensional topology tools needed for it.

\begin{thm}\label{pi1_determines_homeo}
There exist only   finitely many  pairwise  non-diffeomorphic, compact orientable $3$-manifolds without spherical boundary components with given fundamental group $G$.
\end{thm}

\noindent Corollary \ref{fnt_HT} then follows from the fact that the  fundamental groups of   Riemannian manifolds in the class $\mathscr M_{ng}^{\partial}(E,D)$ belong to a finite collection.

  \vspace{5mm}
\subsection{Ramified coverings}\label{ssramifiedcov}${}$ 

\noindent We briefly recall  the construction of  a cyclic  ramified  covering of  a hyperbolic manifold, according  to Gromov-Thurston \cite{GT}.\\
Let $Z_0$ be a two-sided  hypersurface  with boundary in some closed,  orientable n-manifold $X_0$, and call $R_0=\partial Z_0$ the (possibly disconnected) boundary. \\
Cut $X_0$ along $Z_0$, thus obtaining a topological, compact manifold $\dot{X}_0$  with  boundary; the boundary is given by two copies $Z_0^-$, $Z_0^+$ of $Z_0$, with $Z_0^- \cap Z_0^+ =R_0$.
Then, consider the topological manifold $\dot{X}_k$ obtained by taking $k$ copies  $\dot{X}_i$ of $\dot{X}_0$,  for  $i=1,...k-1$, and  gluing  $\dot{X}_i$ to $\dot{X}_{i+1}$ by identifying the boundaries $Z_i^+$, $Z_{i+1}^-$; finally, let $X_k$ be the closed manifold obtained by identifying  $Z_k^+$ to $Z_{1}^-$, and call   $Z_i$ the boundaries so identified inside $X_k$.\\  The resulting manifold $X_k$ can be given a smooth structure with a smooth projection onto the initial manifold $p: X_k \to X_0$ which is a smooth $k$-sheeted covering outside the ramification locus $R=p^{-1}(R_0)$; the ramification locus $R$  is  the boundary of each $Z_i$, and around $R$  the projection  writes  as $(x,z) \mapsto  (x, z^k)$, 
with respect to suitable coordinates for $X_k$ and for $X_0$,   identifying the (trivial) normal bundles of $R$ and $R_0$ to $R_0 \times D^2$. \\
Moreover, choosing  $X_0$ hyperbolic  and the submanifolds $Z_0$, $R_0$ totally geodesic in $X_0$, the new manifold $X_k$ can be given a singular,  locally $CAT(-1)$-metric $g_k$ which makes of  $\bigcup_i Z_i$ a {\em totally geodesic (singular) hypersurface} of $X_k$, with totally geodesic boundary $R$,  and such that  the restriction $p|_{X_k-R}$ is a Riemannian covering  (cp. \cite{GT}). Namely,  $\bigcup_i Z_i$   looks like a $k$-paged book, consisting of $k$ copies of $Z$ joined together at $R$, each pair of consecutive pages forming an angle $2\pi$ and a locally convex subset of $X_k$.
 The singular metric, around the ramification submanifold $R$, can be written  as  \pagebreak
$$\cosh(r) g_0 + k \sinh(r)d\theta^2  +  dr^2$$
 where $r$ represents  the distance to $R$, and $g_0$ the hyperbolic metric of $Z_0$.  \linebreak
(As shown in \cite{GT}, this metric can then be smoothed to obtain a true Riemannian metric $g_k^\epsilon$ of  strictly negative curvature $K(X_k) \leq -1$, and even pinched around $-1$, provided that the normal injectivity radius of $R$ is sufficiently large).\\
We will  call a  Riemannian manifold  obtained by choosing {\em any} Riemannian metric on such $X_k$  (possibly with variable sectional curvature,  of any possible sign, and not necessarily locally isometric to the base hyperbolic manifold $X_0 \setminus R_0$), for any $k \geq 2$, a {\em Riemannian ramified covering  of $X_0$}.
\vspace{4mm}

  {\em Proof of Corollary \ref{corram}.}
	Let $\dot{X}_{k-1} \subset X_k$ be  the union of the  $\dot{X}_i$'s, for $ 2 \leq i \leq k$.\linebreak
	By Van Kampen theorem, $G_k=\pi_1(X_k)$ can be written as the amalgamated product   $G_1\ast_H G_{k-1}$ of $G_1=\pi_1(\dot{X}_1)$ and $G_{k-1}=\pi_1(\dot{X}_{k-1})$ along $H= \pi_1 (Z_1^- \cup Z_1^+)$ (which is immersed in  $G_{k-1}$ via the isomorphisms induced by the identification  diffeomorphisms $Z_1^- \simeq Z_{k}^+$ and $Z_1^+ \simeq Z_2^-$).   
	 As the image $Z_1 \cup Z_k$ of $Z_1^- \cup Z_1^+ $ in $X_k$ is locally convex  with respect to the singular, $CAT(-1)$-metric of $X_k$ described above, the subgroup  $H$ is a malnormal subgroup of $G$ by Proposition \ref{malnormal} of Appendix \ref{appcat(0)}.
%
%
We  conclude that   
  $G$ admits a  $0$-acylindrical splitting; moreover, as $H$ has not index two in $G$ (being malnormal)  the splitting is non-elementary.

\noindent We can then  apply   the  Entropy-Cardinality inequality to a triangular presentation of  $\pi_1(X)$, as in the proof of 
 Theorem \ref{fntn_orb} and \ref{glob_fnt}.
 Namely,  we choose some point $x$ for every $X \in \mathscr R(E,D)$, we  consider a triangular presentation of $\pi_1(X)$ by the  $M$-short generating sets   $S_M$ for $M=2D+\epsilon$,  and as $d(x, g\cdot x) \leq M \cdot |g|_{S_M}$ for all $g \in S_M$, we deduce by Theorem \ref{ent-card_ineq}  that $|S_M |\le N(E,D)$.
Therefore,  the number of possible fundamental groups in  $\mathscr  R(E,D)$ is explictly bounded in terms of $E,D$.  
 Since every $X_k$ is an aspherical manifold admitting a  $CAT(-1)$ metric,  we infer the finiteness of homotopy types (and of diffeomorphisms types in dimension different from $4$) as explained in the proof of Corollary \ref{cormanifoldsk<0}.\qed

\vspace{4mm}
\subsection{Higher dimensional graphs and cusp decomposable manifolds}\label{sscuspdec}${}$ 
\label{sechighergraphs}
%
\vspace{2mm}

{\em Proof of  Corollary \ref{corgraphcusp}}. By  \cite{FLS}, Proposition 6.4, we know that the fundamental groups of irreducible higher graph manifolds admit  $2$-acylindrical splittings. 
On the other hand, the fundamental groups of cusp decomposable manifolds   possess  non-elementary, $1$-acylindrical splittings. Actually, the decomposition of $\pi_1(X)$ corresponding to the cusp decomposition is obtained by identifying the  cusp subgroups, and  these subgroups are malnormal in the fundamental group of each bounded cusp manifold with horoboundary  they belong to, and conjugately separated  if they belong to the same group (cf.  \cite{delHW}, or just apply Proposition \ref{malnormal} of Appendix \ref{appcat(0)} to the whole, convex, cusp neighbourhoods). 
As a consequence, by Lemmas \ref{malnormalAP} and \ref{malnormalHNNext}, the fundamental group of a cusp decomposable manifold  can be presented as a  $0$-step malnormal amalgamated product or HNN-extension, and admits  a  $1$-acylindrical splitting. 
Moreover, the Bass-Serre tree of the splittings corresponding to the  decompositions of  {\em non-elementary} higher graph or cusp decomposable manifolds is, by definition,   neither a vertex nor a line. 
Therefore, the number of fundamental groups of manifolds in the classes   $\mathscr  G ^{\partial}(E,D)$ and   $\mathscr  C (E,D)$  is finite,  by the same argument   used for Theorems \ref{fntn_orb} and  \ref{glob_fnt}.
 Since higher graph and cusp decomposable manifolds are aspherical (cf. \cite{FLS}, Corollary 3.3 and \cite{nguyen}) we immediately infer  the finiteness of the homotopy types  in  $\mathscr  G ^{\partial}(E,D)$ and  $\mathscr  C (E,D)$. \linebreak  By the topological and differential rigidity properties of higher graph and cusp decomposable manifolds recalled in the introduction, we also deduce the finiteness of diffeomorphism types in $\mathscr  G (E,D)$  and   $\mathscr  C (E,D)$.\qed

\begin{rmk}  
	The finiteness result holds, more generally, also for  the diffeomorphism types of  {\em non-irreducible} high dimensional  graph manifold {\em with boundary},  admitting at least one internal walls with transverse fibers, and whose boundary components do not  belong to surface pieces  and  (see  \cite{FLS}, Sections \S5). 
\end{rmk}


\vspace{5mm}

\appendix

\section{Acylindrical splittings of hyperbolic $2$-orbifolds}
\label{apporbifold}

\begin{prop}\label{splitting2orbifolds}
	Let $\mathcal O$ be a compact $2$-orbifold of hyperbolic type with conical singularities. If $\mathcal O$ is not a hyperbolic triangular orbifold, then $\pi_1^{orb}(\mathcal O)$ admits a $2$-acylindrical splitting.
\end{prop}

{\em Proof}.	First notice that all compact  $2$-orbifolds of hyperbolic type with non-empty boundary have orbifold fundamental group which is a  non-trivial  free product of finite and infinite cyclic groups,  hence $\pi_1^{orb}(\mathcal O)$ has a $0$-acylindrical,  splitting in this case. Moreover, the splitting is necessarily non-elementary (otherwise the orbifold would be a disc with two singular points of order two, and it would not have negative orbifold Euler characteristic).\\
Assuming then that $\mathcal O$ is a compact,  $2$-orbifold  of genus $g$ of hyperbolic type without boundary, which is not a hyperbolic triangular orbifold. By the  formula  for the orbifold Euler characteristic one  of the following holds:
	\begin{itemize}
		\item $g>1$;
		\item $g=1$ and $\mathcal O$ has at least one singular point;
		\item $g=0$ and $\mathcal O$ has $m \geq 4$ singular points,   at least one of which has order greater than $2$;
		\item $g=-1$ and $\mathcal O$ has $2$ singular points, one of which of order greater than $2$;
		\item $g=-1$ and $\mathcal O$ has   $m \geq 3$ singular points;
		\item $g=-2$ and $\mathcal O$ has at least one singular point;
		\item $g<-2$.
	\end{itemize}

\noindent The proof then is obtained by cutting any such orbifold   $\mathcal O$ into two $2$-orbifolds with boundary, and using repeatedly the following

\begin{lem}\label{lemmaboundary}
If  $\mathcal O =\mathcal O (g, h;p_1,..., p_k)$ is a compact $2$-orbifold with boundary of hyperbolic type, the infinite cyclic subgroups $\langle d_i \rangle $, corresponding to the boundary curves form a collection of malnormal, conjugately separated subgroups of $\pi_1^{orb} (\mathcal O)$.
\end{lem}

\noindent The lemma can be checked directly by looking at the aforementioned presentations of the orbifold fundamental group: it is sufficient to notice that the boundary curves are represented by primitive elements of infinite order in a non-trivial free product of cyclic groups, different from $\Z_2*\Z_2$,  which do not belong to the same conjugacy class. 
A more geometric justification to manlnormality is that  $\mathcal O $ can be given a geometric structure  of a hyperbolic $2$-orbifold with cusps, with the boundary subgroups $\langle d_i \rangle $ becoming the parabolic  subgroups  associated to the cusps.
\vspace{1mm}

\noindent  Now,   If $g \geq 1$,  choose a simple closed curve $\delta$ which does not disconnect $|\mathcal O|$; after possibly modify the curve $\delta$ in order to avoid the singular points, cut $|\mathcal O|$ along that curve. We obtain a new orbifold $\mathcal O'$ with genus $g-1$ and two new boundary components $\delta_1$, $\delta_2$; clearly,  $\chi_{orb}(\mathcal O')=\chi_{orb}(\mathcal O)<0$. 
   Since $\mathcal O'$ is an orbifold of hyperbolic type, the classes $d_1$, $d_2$, represented by $\delta_1$ and $\delta_2$ in the fundamental group of $\mathcal O'$,  generate  two conjugately separated, malnormal subgroups in $\pi_1^{orb}(\mathcal O')$, by Lemma \ref{splitting2orbifolds}.
Then, by Lemma \ref{malnormalHNNext} we know that  $\pi_1^{orb}(\mathcal O)$ is the HNN-extension $\pi_1^{orb}(\mathcal O')*_{\varphi}$ defined by the isomorphism $\varphi:\langle d_1\rangle\f \langle d_2\rangle$,    $\varphi(d_1)=d_2$; this yields   a  $2$-acylindrical splitting of $\pi_1^{orb}(\mathcal O)$.  By construction, the Bass-Serre tree of this splitting is neither a point nor a line, so the  splitting is   non-elementary.   
   
\noindent    Assume now that  $g=0$ and that $\mathcal O$ has at least $m\geq4$ singular points,   one of which of order $r \geq 3$. Consider a  simple closed curve $\delta$ which separates $|\mathcal O|$ into two disc orbifolds $\mathcal O_1, \mathcal O_2$, each containing at least  $2$ singular points,  with, let's say,   $\mathcal O_1$  containing the singular point of order $r \geq 3$. Denoting by $d$ the classes represented by  the boundary curve  in each $\mathcal O_i$, the  orbifold  fundamental groups have   presentations:
   \small
   $$\pi_1^{orb}(\mathcal O_1)=\langle c_1,..., c_k, d  \; | \; c_1\cdots c_k d =1, \, c_1^{p_1}=\cdots=c_k^{p_k}=1 \rangle$$
   $$ \pi_1^{orb}(\mathcal O_2)=\langle c_{k+1},..., c_m, d   \; | \; c_{k+1}\cdots c_m d =1,  \,  c_{k+1}^{p_{k+1}}=\cdots = c_m^{p_m}=1\rangle$$
   \normalsize
\noindent     Notice that $\chi (\mathcal O_1) = \chi (|\mathcal O_1| ) - 1 - \sum_{i=1}^k (1 - \frac{1}{p_k} )<0$,  therefore the infinite cyclic subgroup $\langle d \rangle$ is a malnormal subgroup of $\pi_1^{orb}(\mathcal O_1)$, by Lemma \ref{splitting2orbifolds}.
  Moreover,   $\pi_1^{orb}(\mathcal O)$ splits non-trivially as $\pi_1^{orb}(\mathcal O_1) \ast_{\langle d \rangle} \pi_1^{orb}(\mathcal O_2) $,  and Lemma \ref{malnormalAP} implies that this  is a non-elementary, $2$-acylindrical splitting.
   
\noindent     If $g=-1$,  consider a closed loop $\delta$ enclosing the singular points $x_1,...,x_m$  and  cut $\mathcal O$ along this loop. 
   We obtain  two orbifolds with boundary: a disc with $m$ singular points $\mathcal O_1$, and a M\"obius strip $\mathcal O_2$.   Observe that $\chi(\mathcal O_1)=\chi (\mathcal O) <0$ and,  calling  again $d$ the classes represented by the boundary loops $\delta$  in each $\pi_1^{orb}(\mathcal O_i)$, the subgroup $\langle d  \rangle $ is   malnormal  in $\pi_1^{orb}(\mathcal O_1)$  by Lemma \ref{lemmaboundary},   while $ \langle d  \rangle$ is a subgroup of   index two in  $\pi_1^{orb}(\mathcal O_2)=\langle a, d \,|\, a^2d =1\rangle$. 
    As 
   $\pi_1^{orb}(\mathcal O)=\pi_1^{orb}(\mathcal O_1) \ast_{\langle d \rangle }\pi_1^{orb}(\mathcal O_2)$, we have  again
 by Lemma \ref{malnormalAP} that $\pi_1^{orb}(\mathcal O)$ is a non-trivial, $1$-step malnormal amalgamated product,  and   possesses a non-elementary, $2$-acylindrical splitting. 
   
  \noindent    Finally,  if $g \leq -2$ then $\mathcal O$ can be cut  along a boundary loop $\delta$ in two orbifolds $\mathcal O_i$, and  we can assume that either  $\mathcal O_1$ has genus $1$ and at least one singular point, or  $\mathcal O_1$ has genus greater than $1$ and no singular points. In the first case, one has $\pi_1^{orb}(\mathcal O_1) = \langle a , c, d  \; | \; a_1^2cd  = c^p=1 \rangle \cong \mathbb{Z} \ast \mathbb{Z}_p $, whereas in the second one $\pi_1^{orb}(\mathcal O_1) = \langle a_1, a_2,   d  \rangle \cong \mathbb{Z}\ast \mathbb{Z}$. In both cases, $\pi_1^{orb}(\mathcal O)$ splits as a non-trivial amalgamated product $\pi_1^{orb}(\mathcal O_1) \ast_{\langle d \rangle }\pi_1^{orb}(\mathcal O_2)$ with $\langle d \rangle$ malnormal in $\pi_1^{orb}(\mathcal O_1) $, which gives again a non-elementary $2$-acylindrical splitting.
\qed

\vspace{2mm}
\section{3-manifolds with prescribed fundamental group}
\label{appJS}

The following statement  is consequence of a number of classical results, which we will recall hereafter for the convenience of the reader:

	\begin{thm} 
		\label{pi1_determines_homeo}
		There exist only  finitely many  pairwise non-diffeomorphic, compact orientable $3$-manifolds without spherical boundary components with given fundamental group $G$.
	\end{thm}

	\vspace{1mm}
	
	To begin with, recall  that  in dimension $3$   the homeomorphism type determines the diffeomorphism type,  by  the celebrated works of Moise, Munkres and Whitehead \cite{Moi}, \cite{Mun}, \cite{Mun2}, \cite{Whi}. 
	
	Now, Theorem \ref{pi1_determines_homeo}  is well-known for {\em closed} $3$-manifolds.
	Actually,  if  $X$ and $X'$ are {\em prime}, closed,  orientable,   
	$3$-manifolds with isomorphic fundamental groups, then $X$ and $X'$ are homeomorphic, unless $X$ and $X'$ are lens spaces; this follows from basic facts of $3$-dimensional topology and from the solution of the Geometrization Conjecture
	(see, for instance, \cite{AFW}, chapters 1\&2).
	Moreover,  by the classification of lens spaces, for every fixed $p\in \N$ there exists only a finite number of lens spaces $L(p,q)$ having $\Z_p$ as fundamental group (see for example \cite{AFW}, pp. 27--28).
	On the other hand, for {\em non-prime}, closed $3$-manifolds, the statement  follows by   Kneser's theorem and the fact that the homeomorphism type of a connected sum is determined by the  prime factors  up to a finite number of choices,   the indeterminacy being  given by the orientations
	of  the summands.
	
	
	The proof of Theorem \ref{pi1_determines_homeo}  for  general compact, orientable $3$-manifolds   with boundary is more tricky and due to Johannson (\cite{jo2}, Corollary 29.3) in the case of {\em irreducible} manifolds {\em with incompressible boundary}\footnote{Johannson's statement is more general and requires the manifolds to be Haken.}. Recall that a compact $3$-manifold $X$ is   {\em irreducible} if any embedded $2$-sphere bounds a $3$-ball. The same result was proved, independently, by Swarup (\cite{swa}, Theorem A), without the incompressibility assumptions. However, the part of Swarup's proof  dealing with possibly compressible boundary components invokes a proposition from \cite{jo3} (namely,  Proposition 3.9) that we were not able to track; since we  noticed that this result, in more recent references like \cite{AFW}, is stated only for irreducible compact $3$-manifolds with incompressible boundary, we find worth filling the details of the proof for general, compact manifolds with boundary without spherical boundary components, assuming Johannson's statement. We will closely follow  Swarup's ideas, so  no claim  of originality is made.
	\vspace{2mm}
	
	Let us recall some basic terminology about $3$-manifolds $X$ with boundary. \\
	A closed, {\em properly embedded} $2$-disk $D\subset X$ (that is, such  that $\partial D \subset \partial X$) is called \textit{essential} if the loop $\partial D$ does not bound any embedded disk in $\partial X$.
	Two such disks $D, D'$ are said to be   {\it parallel} if there is an ambient isotopy sending $D$ into $D'$; 
	and  one says that $X$ has \textit{incompressible boundary} if there are no essential disks.  
	\vspace{2mm}

	\noindent {\em The surgery procedure for irreducible manifolds with compressible boundary.} \\
	Let $X$ be a compact, irreducible $3$-manifold: a {\it disk system} for $X$ is a collection $\mathcal C$ of pairwise disjoint and non-parallel  essential disks;  the system is {\em maximal} if any collection $\mathcal C'$ of essential disks properly containing $\mathcal C$ contains a pair of parallel disks. \linebreak
	Assume that $X$ has compressible boundary: we can then choose a non-empty, maximal disk system  $\mathcal C=\{D_1,...,D_r\}$ and  remove these  disks from $X$. 
	This procedure chops our irreducible manifold $X$ into a finite collection  $\Gamma(X,\mathcal C)$ of irreducible $3$-manifolds with incompressible boundary $X_1,.., X_n$ and finitely many $3$-dimensional balls $B_1,..., B_m$. 
	Moreover, the collection $\Gamma(X, \mathcal C) =\{X_1,.., X_n, B_1,..., B_m\}$ can be given a graph structure:  the  edges $d_i$ of $\Gamma(X,\mathcal C) $ are   in bijection with  the disks $D_i$ of the maximal disk system $\mathcal C$, and two vertices  $v,v' $ of  $\Gamma(X,\mathcal C) $  (possibly with $v=v'$) are connected by  $d_i$  if the  disk $D_i$   bounds the corresponding manifolds or balls.\\
	The irreducible components $X_1,...,X_n$ with incompressible boundary are uniquely determined up to diffeomorphism and do not depend on the particular maximal disc system $\mathcal C$ (see \cite{mat}  pp. 167--168, or   \cite{martelli}); on the other hand, the number $k$ of balls arising from the surgery procedure may depend on the choice of $\mathcal C$. \\
	This procedure can be inverted: we can reconstruct the manifold $X$  from   $\Gamma(X, \mathcal C)$  by gluing back a $1$-handle, i.e. a copy of $D^2\times [0,1]$, for every edge of the graph: roughly speaking, $X$ appears as a ``solid  graph'' whose vertices are the manifolds in the collection $\{X_1,.., X_n, B_1,..., B_m\}$ and whose edges are $1$-handles connecting two (possibly equal) boundary components of the vertices. 
	Using Van Kampen's theorem we see that the fundamental group of $X$ is isomorphic to a free product $\pi_1(X)\cong\pi_1(X_1)*\cdots*\pi_1(X_n)* \mathbb F_k$ where $k$  is the number (possibly zero) of cycles in the graph $\Gamma(X, \mathcal C)$. 
	\vspace{2mm}
	
	Clearly, the number $n$ of  compact, irreducible $3$-manifolds $X_i$ with incompressible boundary components obtained by the surgery procedure is bounded, by Grushko's theorem,  by $N=n+k$. The next Lemma gives a bound of the numbers $r$ and $m$ of, respectively, $1$-handles and balls appearing from the surgery procedure: 
	
	
	\begin{lem}\label{ballsandhandles}
		Let $N$ be the number of irreducible factors of $G = \pi_1(X)$ as a free product: then,  $m \leq 2N$ and $r \leq 3N$.
	\end{lem}

{\em Proof.}		We associate to $\Gamma(X, \mathcal C)$ a graph of groups $\mathscr G (X, \mathcal C)$, by assigning the group  $G_{X_i}=\pi_1(X_i)$   to each vertex $X_i$,  and the trivial groups to the vertices $B_j$ and to every  edge $d_i$. Then, $\pi_1(\mathscr G_Y)\cong  G_{ X_1}*\cdots *G_{ X_n}* \mathbb F_k\cong G$ exactly.\\
		Notice  that, from the  non-parallelism  condition,  the degree of the vertices of $\Gamma(X, \mathcal C)$
		corresponding to the $3$-balls is at least $3$, unless the initial manifold was a solid torus, in which case the graph is just a loop  and the collection of manifolds obtained after the surgery consists of a single $3$-ball; therefore, we may assume that $deg( B_i)\ge 3$ for  $i=1,.., m$. 
		On the other hand,  since the initial manifold has compressible boundary, we know as well that $deg( X_i)\ge 1$ for each $i=1,..,n$. 
		Now, consider a maximal tree $\mathcal T$ in  $\Gamma(X, \mathcal C)$: the maximal tree will have $n+m$ vertices and $n+m-1$ edges. Let $\mathsf E'=\mathsf E(\Gamma(X, \mathcal C))\setminus\mathsf E(\mathcal T)$. Observe that, by construction, the adjunction of each edge of $\mathsf E'$ to $\mathcal T$ corresponds to add a free factor isomorphic to an infinite cyclic group, so $\# \mathsf E' =k$. Then,
		\vspace{-4mm}
		
		\small
		$$2\cdot k+2\cdot (n+m-1)\ge 2\cdot\# \mathsf E(\Gamma(X, \mathcal C))
		=\sum_{i=1}^n deg( X_i)+\sum_{i=1}^{m} deg( B_i)\ge n+3m$$
		\normalsize
		
		\vspace{-2mm}
		\noindent hence  $m\le (2k+n-2)$ and $r=\# \mathsf E(\Gamma(X, \mathcal C))\le k+n+m-1\le 3k+2n-3$ 
		and we conclude that  $m$ and $r$ are (roughly) bounded respectively by  $2N$ and $3N$.
\qed
	
	\vspace{3mm}
 {\em Proof of Theorem \ref{pi1_determines_homeo} for orientable manifolds with boundary}. 
	Let $G$
	be a (compact) $3$-manifold group, whose decomposition as a free product has $N$ indecomposable factors.
	If $X$ is an orientable,  compact manifold with boundary without spherical boundary components and fundamental group $G$, it has a prime decomposition as a connected sum   of irreducible manifolds and copies of $S^2 \times S^1$ (the only prime, non irreducible manifold without spherical boundary components), with at most  $N$ factors. The homeomorphism type of a connected sum being uniquely determined by its factors and their orientations,  it will then be enough to prove the theorem for irreducible manifolds.
	Now, by Lemma \ref{ballsandhandles}, any compact, irreducible $3$-manifold $X$ with fundamental group isomorphic to $G$ can be splitted using the surgery procedure in   $n\leq N$ irreducible  $3$-manifolds $X_i$ with incompressible boundary and   fundamental group $G_i$, plus a number  $m \leq 2N$ of  $3$-balls; and $X$ is obtained as a solid graph on these pieces,   attaching  at most $r \leq 3N$ 1-handles.
	Moreover, notice that, by Kneser's Theorem (holding for irreducible $3$-manifolds with incompressible boundary components), the fundamental group of each $X_i$ is indecomposable, hence isomorphic to one indecomposable factor of the free product decomposition of $G$.
	Now, by Johannson's theorem, for each indecomposable factor $G_i$ of $G$ 
	there exist only finitely many  non-homeomorphic  irreducible  $3$-manifolds $X_{i,\alpha}$ with incompressible boundary with fundamental group $G_i$.
	Moreover, any two disks $D, D'$ in one of these $X_{i, \alpha}$ are isotopic,
and there are only two isotopy classes of diffeomorphisms $D^2\f D^2$ (corresponding to the identity and to a reflection with respect to one axis); hence, once fixed two such pieces  $X_{i, \alpha}$ and $X_{j, \beta}$, there are essentially two inequivalent ways of attaching a 1-handle to them. Therefore, there are only finitely many manifolds which can be obtained as a solid graph on the (finitely many) pieces $X_{i, \alpha}$, which concludes the proof.\qed

 \vspace{1mm}
\section{Malnormal subgroups of  CAT(0)-groups}
\label{appcat(0)}

We start giving a method to detect malnormal subgroups in fundamental group of locally $CAT(0)$-spaces. 

\begin{prop} \label{malnormal}
	Let $Z$ be a  compact,   locally convex subspace of a compact, complete locally  $CAT(0)$-space $X$.
	Assume that $X$ is negatively curved around $Z$: then, $H=\pi_1 (Z)$ is malnormal in   $G=\pi_1 (X)$.
\end{prop}

By {\em negatively curved around $Z$} we mean that $Z$ has a neighbourhood  $U(Z)$ in $X$  such that $U(Z)\setminus Z$   is a locally $CAT(-k)$-space, for some $k>0$.
Notice that this covers the case where  $X$ is a complete Riemannian with sectional curvature $k_X < 0$, with no a-priori negative upper bound on the curvature.\\

{\em Proof}.
	Let  $ \tilde X \rightarrow X$ the universal covering map: $\tilde X$ it is a $CAT(0)$-space. \linebreak
	Let $\tilde z_0 \in p^{-1}(Z) \subset \tilde X$ be a point projecting to $z_0 \in Z$, and $H=\pi_1(Z,z_0)$. 
	Finally, let  $C_{\tilde z} p^{-1}(Z)$   denote the connected component of $p^{-1}(Z)$ containing the point $\tilde z$, and    
	$\tilde Z_{\tilde z_0}$ the subset of $ \tilde X$ obtained by lifting from $\tilde z_0$ any curve $\gamma$ of $Z$ based at $z_0$, and taking the endpoint $\widetilde{\gamma} (1)$ of the lift $\widetilde \gamma$. Then:
	\vspace{1mm}
	
	(1)\, {\em  $C_{\tilde z_0} p^{-1}(Z) = \tilde Z_{\tilde z_0}$, and is a covering  of $Z$.} \\
	The fact that $C_{\tilde z} p^{-1}(Z)$ is a covering follows from ordinary theory of coverings, and the inclusion $\tilde Z_{\tilde z_0} \subset C_{\tilde z_0} p^{-1}(Z)$ is trivial. 
	Conversely, if $\tilde z  \in C_{\tilde z_0} p^{-1}(Z)$, it can be joined to $\tilde z_0$ by a curve $\tilde \gamma$ whose projection $\gamma$ stays in $Z$; hence $\tilde z  = \tilde \gamma (1) \in \tilde Z_{\tilde z_0} $ by definition.
	\vspace{1mm}
	

	(2)\, {\em  $\tilde Z_{\tilde z_0}$  is the universal covering of $Z$, and $H$ injects in $G$. }  \\
 Actually, since locally $CAT(0)$ spaces are locally convex
\footnote{$CAT(0)$ spaces are assumed to be locally geodesic spaces (though non necessarily geodesic spaces), by definition.}
, every class in $\pi_1(X,z_0)$ can be realized by a locally geodesic loop. Now,   every locally geodesic  loop $\gamma$   representing a class of $\pi_1(Z,z_0)$ lifts to a local geodesic  $\tilde \gamma$ of $\tilde X$ from $\tilde z_0$ (the covering map being locally isometric). But every local geodesic in  a $CAT(0)$ space is a true geodesic, hence $\tilde \gamma$  is not closed: this shows that $\tilde Z_{\tilde z_0}$ is simply connected, and that $\gamma$ does not represent the trivial element of $G$.
	\vspace{1mm}  
	
	(3)\, {\em  $\tilde Z_{\tilde z_0}$, endowed with the length structure induced by $Z$,   is  isometrically embedded in $\tilde X$; therefore, it is a convex subset of $\tilde X$.} \\
  In fact, since $Z$ is locally convex in $X$,   the inclusion $\tilde Z_{\tilde z_0} \subset \tilde X$  is a local isometry;  but, $\tilde  X$ being $CAT(0)$, geodesics in  $\tilde  X$ are unique, which implies that $\tilde Z_{\tilde z_0}$ is convex in $\tilde X$ and that the inclusion is a true isometric embedding.
 \vspace{1mm} 
		

	(4)\, {\em   $C_{g \tilde z_0} p^{-1}(Z)  = g \cdot \tilde Z_{\tilde z_0} $.}\\
	As in (1) one sees that $C_{g \tilde z_0} p^{-1}(Z) =  \tilde Z_{g \tilde z_0} $ (the subset obtained by lifting from $g \tilde z_0$ any curve $\gamma$ of $Z$ with base point $z_0$), which clearly equals  $g \cdot \tilde Z_{g  \tilde z_0}$.
	\vspace{1mm}
	
	(5)\, {\em    $\Stab_G (\tilde Z_{\tilde z_0}) = H$.} \\
	The elements of $H$ clearly stabilize $\tilde Z_{\tilde z_0}$ (recall that $h \in H$ acts on $\tilde x \in \tilde X$ by lifting from $\tilde z_0$ the composition of a geodesic $c$ from $z_0$ to $x=p(\tilde x)$ with a loop $\gamma$ at $z_0$ representing $h$; so, the final point of the lift $\widetilde{ \gamma c }$ belongs to $\tilde Z_{\tilde z_0}$ by definition of $\tilde Z_{\tilde z_0}$).
	Conversely: if $g \in \Stab_G (\tilde Z_{\tilde z_0})$, then $g \tilde z_0 \in \tilde Z_{\tilde z_0}$, and then the geodesic $\tilde \gamma$ joining $\tilde z_0$ to $g \tilde z_0$ stays in $\tilde Z_{\tilde z_0}$ (since this is a convex subset of $\tilde X$). As $g$ is represented by the projection $\gamma$ of $\tilde \gamma$  in $X$, which is included in $Z$, then $g \in \pi_1(Z)=H$.
	\vspace{1mm}
	
	(6)\, {\em   $\Stab_G (g \cdot \tilde Z_{\tilde z_0}) = gHg^{-1}$, and  the number of connected components of  $p^{-1}(Z)$ is in bijection with the cosets space $ G / H$.}\\
	Both assertions follow  from (4) and (5).
	\vspace{1mm}
	
	(7)\,  {\em  Every  $h \in H$   acts  on $\tilde X$ by hyperbolic isometries, and  the subset  $Min(h)$ where the displacement function $d(\tilde x, h \tilde x)$ attains its minimum 
	 is   included in $\tilde Z_{\tilde z_0}$.}\\
	Since  the action of $G=\pi_1 (X, z_0)$ on $\tilde X$ is cocompact and without fixed points,  then every element  of $G$ acts on $\tilde X$ by hyperbolic isometries.
	We shall now prove that $Min(h)$ is entirely included in $ \tilde Z_{\tilde z_0}$,  for every  $h \in H$. 
	Actually, let $\tilde x_0$ be an arbitrary point of minimum for the displacement function $s_h( \tilde x) = d(\tilde x, h \tilde x)$, and consider the projection $p: \tilde X \to \tilde Z_{\tilde z_0}$ (this  is well defined, since $\tilde Z_{\tilde z_0}$ is a convex subset). As $\tilde Z_{\tilde z_0}$  is invariant under $h$ by (5), and  since $p$ is a projection,  we have $h \cdot p(\tilde x_0)= p(h \cdot \tilde x_0)$. Therefore
	$$d(p(\tilde x_0),h \cdot p(\tilde x_0)) = d(p(\tilde x_0),  p(h \cdot \tilde x_0)) \leq  d(\tilde x_0, h \cdot \tilde x_0)$$
	This shows that the point  $p(\tilde x_0)\in \tilde Z_{\tilde z_0}$ also realizes the minimum of $s_h(x)$.
	By  the $h$-invariance and the convexity of  $Min(h)$, we deduce that the orbits $\{ h^n \cdot \tilde x_0 \}$ and $\{ h^n \cdot \tilde p(x_0) \}$ lie on two parallel geodesics $\gamma$ and $p(\gamma)$, entirely included in $Min(h)$; moreover, $p(\gamma) \subset  \tilde Z_{\tilde z_0}$, as $ \tilde Z_{\tilde z_0}$ is convex.
	So, $Min(h)$ contains a flat band bounding $\gamma$ and  $p(\gamma)$;  the lifted neighbourhood $\tilde U(\tilde Z_{\tilde z_0})$ of $\tilde Z_{\tilde z_0}$  being strictly negative curved outside $\tilde Z_{\tilde z_0}$, this shows that $\gamma$ and  $\tilde x_0$  are necessarily included in $\tilde Z_{\tilde z_0}$.

	\vspace{1mm} 
	

	(8)\, {\em  $H$ is malnormal in $G$.} \\
	Assume that there exists   $h \in H^\ast$ and  $g \in G$ such that $h'=ghg^{-1} \in H$. \linebreak
	By (7), $Min(h)$ is included in $  \tilde Z_{\tilde z_0}  = C_{\tilde z_0} p^{-1}(Z)$; but as $h'=ghg^{-1}$ is in $H$, we also have  $Min(h') \subset   C_{\tilde z_0} p^{-1}(Z)$. However,  $Min(h')=Min(ghg^{-1})=g\cdot Min(h)$ is  included in $ g\cdot \tilde Z_{\tilde z_0} =   C_{g \tilde  z_0} p^{-1}(Z) $, which is disjoint from $C_{\tilde z_0} p^{-1}(Z)$ if $g \not\in H$, by (6). \linebreak
	This shows that $g \in H$ and that $H$ is malnormal in $G$.
\qed

	\normalsize
\end{document}